\documentclass[a4paper,11pt]{article}
\usepackage{amssymb}
\usepackage{amsmath}
\usepackage{amsthm}
\usepackage{graphicx}
\newtheorem{thm}{Theorem}[section]

\newcommand{\MLONE}{E_{\alpha,1}}
\newcommand{\Hgam}{H_{\alpha}(0,T;L^2(\OOO))}
\newcommand{\Ahalf}{A_0^{\frac{1}{2}}}
\newtheorem{proposition}[thm]{Proposition}
\newtheorem{remark}[thm]{Remark}
\newtheorem{example}[thm]{Example}
\newtheorem{lemma}[thm]{Lemma}
\newtheorem{definition}[thm]{Definition}
\newtheorem{theorem}[thm]{Theorem}

\numberwithin{equation}{section}

\newcommand{\sumn}{\sum_{n=1}^{\infty}}
\newcommand{\QQQQ}{\Omega \times (0,T)}
\newcommand{\UPP}{\overline{u}}
\newcommand{\ep}{\varepsilon}
\newcommand{\la}{\lambda}
\newcommand{\va}{\varphi}
\newcommand{\ppp}{\partial}

\newcommand{\DDDA}{\DDD(A_0^{\gamma})}

\newcommand{\pppa}{\partial_t^{\alpha}}
\newcommand{\OOO}{\Omega}
\newcommand{\sumij}{\sum_{i,j=1}^d}
\newcommand{\NUNU}{\ppp_{\nu_A}}
\newcommand{\ooo}{\overline}
\newcommand{\DDD}{\mathcal{D}}
\newcommand{\hhalf}{\frac{1}{2}}
\newcommand{\ddda}{d_t^{\alpha}}

\newcommand{\UPPY}{\ooo{y}}
\newcommand{\LOWY}{\underline{y}}
\newcommand{\lowa}{\underline{a}}
\newcommand{\LOP}{\underline{u}}
\newcommand{\R}{\mathbb{R}}
\newcommand{\C}{\mathbb{C}} 
\newcommand{\N}{\mathbb{N}}

\newcommand{\CC}{_{0}C^{\infty}(0,T)}
\newcommand{\HH}{H_{\alpha}}

\allowdisplaybreaks

\title{Comparison principles for the time-fractional diffusion equations 
with the Robin boundary conditions. \\ Part II: Semilinear equations}
\author{
Yuri Luchko \thanks{Department of Mathematics, Physics, and Chemistry, Berlin University of Applied Sciences and Technology,   Luxemburger Str. 10, 
13353 Berlin, Germany,
e-mail: {\tt luchko@bht-berlin.de}} \space and  \space
 Masahiro Yamamoto
\thanks{ Department of Mathematical Sciences, 
The University of Tokyo,
Komaba, Meguro, 
153-8914 Tokyo, Japan and 
Department of Mathematics, Faculty of Science, Zonguldak B\"ulent Ecevit 
University, 
Zonguldak 67100, T\"urkiye, 
e-mail:
{\tt myama@ms.u-tokyo.ac.jp}}
}
\date{}
\begin{document}
\maketitle
\begin{abstract}
In this paper, we deal with analysis of the initial-boundary value problems 
for the semilinear time-fractional diffusion equations, while 
the case of the linear equations was considered in the first part of 
the present work. 
These equations contain
uniformly elliptic spatial differential operators of the second order and
the Caputo type fractional derivative acting in the fractional Sobolev 
spaces as well as a semilinear term that depends on the spatial variable, 
the unknown function and its gradient. 
The boundary conditions are formulated in form of the homogeneous 
Neumann or Robin conditions. For these problems, we first prove 
uniqueness and existence of their 
solutions. Under some suitable conditions, we then show the non-negativity 
of the solutions and derive several comparison principles. We also apply 
the monotonicity method by upper and lower solutions to deduce some 
a priori estimates for solutions to the initial-boundary value problems 
for the semilinear time-fractional diffusion equations.
Finally, we consider some initial-boundary value problems for systems of the linear and semilinear time-fractional diffusion equations and prove non-negativity of their solutions under the suitable conditions.
\end{abstract}




\section{Introduction} 
\label{sec:1}

\setcounter{section}{1} 
\setcounter{equation}{0} 

In this paper, we mainly treat the initial-boundary value problems for a 
semilinear time-fractional 
diffusion equation
$$
\pppa (u(x,t)-a(x)) = \sum_{i,j=1}^d \ppp_i(a_{ij}(x)\ppp_j u(x,t))
$$
\begin{equation}
\label{(1.5)}
+ \sum_{j=1}^d b_j(x,t)\ppp_ju(x,t) + c(x,t)u(x,t) 
+ f(u)(x,t),\  x \in \Omega,\, 0<t<T,
\end{equation}                               
where $\pppa$ is the Caputo fractional derivative of order 
$\alpha\in (0,1)$ defined  on the fractional Sobolev spaces 
(see Section \ref{sec2} for details), and $\OOO \subset \R^d$ is
a bounded domain with a smooth boundary $\ppp\OOO$.

In the equation \eqref{(1.5)}, $f$ is a semilinear term, and
$f: V\subset L^2(\OOO)\, \longrightarrow \, L^2(\OOO)$
is a mapping, where $V$ is suitably prescribed.
In what follows, we often use the notation $f(v) = f(v)(x)$ omitting the independent variable $x$.
Two typical examples of the semilinear term $f$ are $f(v) = v(x)^2$ with 
$V = L^{\infty}(\OOO)$ and $f(v) = \sum_{k=1}^d \mu(x)v(x)\ppp_kv(x)$ 
with a fixed $\mu\in L^{\infty}(\OOO)$ and with
$V = H^1(\OOO) \cap L^{\infty}(\OOO)$.
Please note that in these  examples, $f(v)$ means substitution of a function
$v$ into the semilinear term $f(\cdot)$. However, 
 we also cover the case of the semilinear terms $f$ depending on several variables, see  Section 3 for details. In applications, the semilinear term $f$   is responsible for description of the reaction component of the anomalous transport processes modeled by the equation (1.1).
 
In this paper, for the sake of simplicity, we restrict ourselves to the case of the spatial dimensions
$d=1,2,3$. The case $d\ge 4$ requires more advanced arguments and will be not considered here. 

In what follows, we always assume that all functions under consideration 
are real-valued or real vector-valued and the following conditions are 
satisfied:
\begin{equation}
\label{(1.2)}
\left\{ \begin{array}{rl}
& a_{ij} = a_{ji} \in C^1(\ooo{\OOO}), \quad 1\le i,j \le d, \\
& b_j,\, c \in C^1([0,T]; C^1(\ooo{\OOO})) \cap C([0,T];C^2(\ooo{\OOO})),
  \quad 1\le j \le d, \\
& \mbox{and there exists a constant $\kappa>0$ such that}\\
& \sumij a_{ij}(x)\xi_i\xi_j \ge \kappa \sum_{j=1}^d \xi_j^2, \quad
x\in \OOO, \, \xi_1, ..., \xi_d \in \R.
\end{array}\right.
\end{equation}

To formulate the boundary conditions, we define the co-normal derivative
$\NUNU w$ with respect to the operator $\sumij \ppp_j(a_{ij}\ppp_i)$ by 
\begin{equation}
\label{(1.3)}
\NUNU w(x) =  \sumij a_{ij}(x)\ppp_jw(x)\nu_i(x), \quad x\in \ppp\OOO,
\end{equation}
where $\ppp_j = \frac{\ppp}{\ppp x_j}$, $j=1, 2, ..., d$ and 
$\nu = \nu(x) =: (\nu_1(x), ..., \nu_d(x))$ is the unit outward normal
vector to the surface $\ppp\OOO$ at the point $x := (x_1,..., x_d) \in 
\ppp\OOO$. 

For the equation \eqref{(1.5)},  we consider the initial-boundary value 
problems with the initial conditions formulated as some inclusions 
in the fractional Sobolev spaces (see Section \ref{sec2}) and either the 
Neumann boundary condition
\begin{equation}
\label{(1.3a)}
\ppp_{\nu_A}u = 0 \quad \mbox{on $\ppp\OOO \times (0,T)$}   
\end{equation}
or the Robin boundary condition
\begin{equation}
\label{(1.4)}
\ppp_{\nu_A}u + \sigma(x)u = 0 \quad \mbox{on $\ppp\OOO \times (0,T)$}.   
\end{equation}
Throughout this paper, we assume that
$$
\sigma \in C^{\infty}(\ppp\OOO), \quad
\sigma(x) \ge 0 \quad \mbox{for $x\in \ppp\OOO$},
$$
although this  regularity condition can be weakened.
 
The main results presented in this paper are the comparison principles 
and the monotonicity method for the semilinear time-fractional 
diffusion equation \eqref{(1.5)} with the 
boundary condition \eqref{(1.3a)} or \eqref{(1.4)}  and a suitably formulated 
initial condition.  Moreover, we also derive some results for the systems of the semilinear 
time-fractional diffusion equations.

The case of the initial-boundary value problems for the linear time-fractional 
diffusion equations of type \eqref{(1.5)} with the semilinear term $f$ 
depending just on the spatial and time variables, but not on the unknown 
function $u$ has been treated in \cite{LYP1}. In particular, 
under some suitable conditions, the non-negativity of solutions and 
the comparison 
principles for solutions to these problems were deduced there. 

The non-negativity properties, the maximum and comparison principles, and 
their applications are well-known and intensively used for  both linear and 
semilinear  partial differential equations of parabolic type ($\alpha=1$
in the equation \eqref{(1.5)}), see, e.g., \cite{E}, 
 \cite{Fr}, \cite{Pao2}, \cite{PW}, \cite{RR}, and Section 52 in 
\cite{QS}. Here we do not intend to provide comprehensive references.
Moreover, we refer to  \cite{Am},  \cite{Ke},  and \cite{Pao2}
for discussions of the monotonicity method by upper and lower solutions for 
studying properties of the solutions to the initial-boundary value problems 
for the semilinear  partial differential equations of parabolic type. In its 
turn, the monotonicity method is based on a corresponding comparison principle 
and thus, one can apply this technique as soon as one has a comparison 
principle at one's fingertips.   

In \cite{Pi1} and \cite{Pi2}, Pierre  discussed unique existence of solutions to the initial-boundary value problems for the systems of the semilinear parabolic PDEs of type \eqref{(1.5)} with $\alpha=1$ as well as their non-negativity.

However, to the best of authors' knowledge, no results of this kind  for the 
initial-boundary value problems for 
the semilinear time-fractional diffusion equations of type \eqref{(1.5)}  
with the Neumann or Robin 
boundary conditions  were yet presented in the literature. 
The main objective 
of this paper is to contribute to
the theory of the semilinear time-fractional diffusion equations 
and to discuss the non-negativity property and the comparison principles 
for the solutions to the initial-boundary value problems for the equation 
\eqref{(1.5)} with the boundary condition
\eqref{(1.3a)} or \eqref{(1.4)} and with an appropriately formulated initial 
condition (see Section \ref{sec2}).  
The arguments employed for derivation of these results
rely on an operator theoretical approach to the fractional integrals and 
derivatives in the fractional Sobolev spaces which is an extension of the 
theory well-known for the PDEs of parabolic type 
(see, e.g., \cite{He}, \cite{Pa}, \cite{Ta}).

The rest of this paper is organized as follows. In Section \ref{sec2}, 
we present some notations, definitions, and known results needed for our 
analysis in the next sections.  
Section \ref{sec5} is devoted to the results regarding the unique existence 
of solutions to the 
initial-boundary value problems for the semilinear time-fractional  
diffusion equations of type \eqref{(1.5)}. Section \ref{sec6} presents 
the comparison 
principles and some of their applications for derivation of important 
properties of the solutions to the initial-boundary value 
problems for the semilinear time-fractional  diffusion equations.
In Sections 5 and 6, under some suitable conditions, we prove the non-negativity of solutions to the initial-boundary value problems for the systems of the linear 
and semilinear time-fractional diffusion equations.
Finally, in the last section, some conclusions, remarks, and directions 
for further research are formulated.

\section{Preliminaries } 
\label{sec2}

\setcounter{section}{2}
\setcounter{equation}{0}
In this section, we present some notions, notations, and results needed for 
further discussions. 

In order to shorten the formulation of the equation \eqref{(1.5)}, 
for $x \in \Omega, \thinspace 0<t<T$,  we define the operator
\begin{equation}
\label{(2.1)}
-Au(x,t) = \sum_{i,j=1}^d \ppp_i(a_{ij}(x)\ppp_jv(x,t))
+ \sum_{j=1}^d b_j(x,t)\ppp_ju(x,t) + c(x,t)u(x,t),            
\end{equation}
where the coefficients 
$a_{ij}(x), b_j(x,t), c(x,t)$ are assumed to satisfy 
the conditions \eqref{(1.2)}. 
 
In this paper, we deal with the following initial-boundary value problem 
for the semilinear  time-fractional diffusion equation with the Caputo 
fractional derivative 
of order $\alpha\in (0,1)$:
\begin{equation}
\label{(2.2)}
\left\{ \begin{array}{rl}
& \pppa (u(x,t)-a(x)) + Au(x,t) = f(u)(x,t), 
\quad x \in \Omega, \thinspace 0<t<T, \\
& \NUNU u + \sigma(x)u(x,t) = 0, \quad x\in \ppp\OOO, \, 0<t<T,
\end{array}\right.
\end{equation}
and the initial condition \eqref{incon} formulated below. 

In the equation \eqref{(2.2)}, the Caputo fractional derivative of order $\alpha\in (0,\, 1)$ defined by
$$
\ddda w(t) = \frac{1}{\Gamma(1-\alpha)}\int^t_0
(t-s)^{-\alpha}\frac{dw}{ds}(s) ds
$$
is interpreted as an operator  $\pppa$ acting on the fractional Sobolev 
spaces.  For the readers convenience, based on $d_t^{\alpha}$,
we provide below a definition of this operator. 

Henceforth, by $L^2(D)$ and $H^k(D)$ with $k\in \mathbb{N}$,
we mean the usual Lebesgue space and the Sobolev spaces 
on an interval $D\subset (0,T)$ or $D = \Omega$
(\cite{Ad}, \cite{Gri}). 

Furthermore, for a Banach space $X$,
we consider a function $t\in (0,T) \, \longrightarrow\,u(t) \in X$, see, e.g., pp. 351--352 in \cite{E} or 
\S2 in Chapter 1 of \cite{Yagi}.
Then, it is well-known that $L^2(0,T;X)$ is  a Banach space provided that 
$X$ is a Banach space.

Now, for $\alpha>0$, we define
the Riemann-Liouville fractional integral operator $J^{\alpha}$ 
for an $X$-valued function $u \in L^2(0,T;X)$ by 
$$
J^{\alpha}u(t) 
:= \frac{1}{\Gamma(\alpha)}
\int^t_0 (t-s)^{\alpha-1} u(s) ds, 
$$
and we regard $J^{\alpha}u$ as a function acting from $(0,T)$ to $X$.  
We set $J^0 = I$, where $I$ is the identity operator on $X$.

Then the following statement holds valid (see, e.g., \cite{HY} for the proof): 

\begin{proposition}
\label{1}

Let $\alpha>0$. Then

\noindent 
$\mathrm{(i)}$ The operator $J^{\alpha}: L^2(0,T;X) 
\longrightarrow L^2(0,T;X)$ is bounded and 
injective.

\noindent $\mathrm{(ii)}$ $J^{\alpha+\beta} = J^{\alpha}J^{\beta}$ 
in $L^2(0,T;X)$ for $\alpha, \beta>0$.
\end{proposition}

Therefore, the inverse operator to $J^{\alpha}$, denoted by $J^{-\alpha}$,
exists and it acts from $J^{\alpha}L^2(0,T;X)$ to 
$L^2(0,T;X)$. 

Then we define the space
$$
\left\{ 
\begin{array}{rl}
& H_{\alpha}(0,T;X) := J^{\alpha}L^2(0,T;X), \\
& \Vert v\Vert_{H_{\alpha}(0,T;X)}:= \Vert J^{-\alpha	}v\Vert
_{L^2(0,T;X)}
\quad \mbox{for } \ v\in H_\alpha(0,T;X).
\end{array}
\right.
$$
In other words, the equality  $\Vert J^{\alpha} u\Vert_{H_{\alpha}(0,T;X)}
:= \Vert u\Vert_{L^2(0,T;X)}$ holds true for $u\in L^2(0,T;X)$. 
Moreover, the space 
$H_{\alpha}(0,T;X)$ is a Banach space with the norm $\Vert \cdot\Vert
_{H_{\alpha}(0,T;X)}$
(see, e.g., \cite{HY} for the proof).

Now we define the fractional derivative $\pppa$ as follows:
$$
\ppp_t^{\alpha}:= (J^{\alpha})^{-1}, \quad \DDD(\pppa) 
= H_{\alpha}(0,T;X).
$$
By the definition of $H_\alpha(0,T;X)$, we can 
readily see that $J^{\alpha}: L^2(0,T;X) \longrightarrow 
H_{\alpha}(0,T;X)$ 
is injective and surjective, and 
$$
\pppa J^{\alpha}u = u \quad \mbox{for all $u \in L^2(0,T;X)$}, \quad
J^{\alpha}\pppa v = v \quad \mbox{for all $v \in H_{\alpha}(0,T;X)$ }
$$
and
$$
\Vert \pppa u\Vert_{L^2(0,T;X)} \sim \Vert u\Vert_{H_{\alpha}(0,T;X)}
\quad \mbox{for all $u \in H_{\alpha}(0,T;X)$}.
$$
Here and henceforth, we denote the norm equivalence by $\sim$. 
It is worth mentioning that $\frac{d}{dt}$ means the 
pointwise differentiation and $\ppp_t^1$ is an operator with the 
domain $H_1(0,T;X)$, where the operation is the same as 
$\frac{d}{dt}$.

For the reader's convenience, we formulate some of the properties of the operators introduced above needed for the further discussions (see, e.g., \cite{HY} for the proofs). 

\begin{proposition}
\label{2}
Let $\alpha>0$, $\beta\ge 0$. Then

\noindent $\mathrm{(i)}$ $J^\alpha:$ $H_\beta(0,T;X)\longrightarrow 
H_{\alpha+\beta}(0,T;X)$ 
is injective and surjective, and \\
$\|J^\alpha u\|_{H_{\alpha+\beta}(0,T;X)} \sim \|u\|_{H_\beta(0,T;X)}$ for 
$u\in H_\beta(0,T;X)$. 

\noindent $\mathrm{(ii)}$ $\partial_t^\alpha:$ $H_{\alpha+\beta}(0,T;X)
\longrightarrow 
H_{\beta}(0,T;X)$ is injective and surjective, and \\
$\|\partial_t^\alpha u\|_{H_{\beta}(0,T;X)} 
\sim \|u\|_{H_{\alpha+\beta}(0,T;X)}$ for $u\in H_{\alpha+\beta}(0,T;X)$. 
\end{proposition}
\begin{proposition}
\label{3}

$\mathrm{(i)}$ Let $\alpha,\beta>0$. 
Then $\partial_t^{\alpha+\beta} u = \partial_t^{\alpha}\partial_t^{\beta} u 
= \partial_t^{\beta}\partial_t^{\alpha} u$ for 
$u\in H_{\alpha+\beta}(0,T;X)$.

\noindent $\mathrm{(ii)}$ Let $0 < \beta \le \alpha$. Then 
$J^{\alpha-\beta} u = \partial_t^\beta J^\alpha u$ for 
$u\in L^2(0,T;X)$ and \\
$J^{\alpha-\beta} u = J^\alpha \partial_t^\beta u$ for 
$u\in H_\beta(0,T;X)$.

\noindent $\mathrm{(iii)}$ Let $0<\alpha<\beta$. Then 
$\partial_t^{\beta-\alpha} u = \partial_t^\beta J^\alpha u$ for
$u\in H_{\beta-\alpha}(0,T;X)$ and \\
$\partial_t^{\beta-\alpha} u  = J^\alpha \partial_t^\beta u$ for 
$u\in H_\beta(0,T;X)$.
\end{proposition}

In the case of $X=\R$,  a characterization of the space 
$H_{\alpha}(0,T;\R)$ in terms of the Sobolev-Slobodeckij spaces was provided in \cite{GLY} and \cite{KRY} and we represent it in what follows. 

We start with the definition of the Sobolev-Slobodeckij space 
$H^{\alpha}(0,T)$ by 
$H^{\alpha}(0,T) = \{ v;\, \Vert v\Vert_{H^{\alpha}(0,T)} < \infty\}$,
equipped with the norm
$$
\Vert v\Vert_{H^{\alpha}(0,T)}:=
\left( \Vert v\Vert^2_{L^2(0,T)}
+ \int^T_0\int^T_0 \frac{\vert v(t)-v(s)\vert^2}{\vert t-s\vert^{1+2\alpha}}
dtds \right)^{\hhalf}.
$$
Then we have the relation (see, e.g., \cite{GLY}, \cite{KRY})
$$
H_{\alpha}(0,T) = 
\left\{ \begin{array}{rl}
&H^{\alpha}(0,T), \quad  0<\alpha<\hhalf, \\
&\left\{ v \in H^{\hhalf}(0,T);\, \int^T_0 \frac{\vert v(t)\vert^2}{t}
dt < \infty \right\}, \quad  \alpha=\hhalf, \\
& \{ v \in H^{\alpha}(0,T);\, v(0) = 0\}, \quad \hhalf < \alpha < 1
\end{array}\right.
$$
with the norm defined by
$$
\Vert v\Vert_{H_{\alpha}(0,T)} = 
\left\{ \begin{array}{rl}
&\Vert v\Vert_{H^{\alpha}(0,T)}, \quad  \alpha \ne \hhalf, \\
&\left( \Vert v\Vert_{H^{\hhalf}(0,T)}^2
+ \int^T_0 \frac{\vert v(t)\vert^2}{t}dt\right)^{\hhalf}, \quad
\alpha=\hhalf
\end{array}\right. .
$$

Moreover, the operator $\pppa$ is the minimum closed extension of 
$\ddda$ with the domain $\DDD(\ddda) = \CC: = 
\{ u \in C^1[0,T];\thinspace u(0) = 0\}$ and the representation 
$$
H_{\alpha}(0,T):= \ooo{\CC}^{H_{\alpha}(0,T)}
$$ 
holds true (\cite{KRY}).  

Having defined the Caputo fractional  derivative on the space 
$H_{\alpha}(0,T;L^2(\OOO))$, we can introduce a suitable 
initial condition for the problem 
\eqref{(2.2)} in form of the following inclusion:
\begin{equation}
\label{incon}
u - a \in \Hgam.
\end{equation}

The complete formulation of the initial-boundary value problem 
for the semilinear  time-fractional diffusion equation \eqref{(1.5)} takes now 
the following form:  
\begin{equation}
\label{(2.3)}
\left\{ \begin{array}{rl}
& \pppa (u(x,t)-a(x)) + Au(x,t) = f(u)(x,t), 
\quad x \in \Omega, \thinspace 0<t<T, \\
& \NUNU u(x,t) + \sigma(x)u(x,t) = 0, \quad x\in \ppp\OOO, \, 0<t<T,\\
& u - a \in H_{\alpha}(0,T;L^2(\OOO)).
\end{array}\right.
\end{equation}

It is worth mentioning that the term $\pppa (u - a)$  in the first line of 
\eqref{(2.3)} is 
well-defined due to inclusion formulated in the third line of \eqref{(2.3)}.
In particular, for $\frac{1}{2} < \alpha < 1$, the Sobolev embedding leads 
to the inclusions
$\HH(0,T;L^2(\OOO)) \subset H^{\alpha}(0,T;L^2(\OOO))\subset 
C([0,T];L^2(\OOO))$. 
This means that $u\in \HH(0,T;L^2(\OOO))$ implies  $u \in C([0,T];L^2(\OOO))$  
and thus in this case the initial condition can be formulated as 
$u(\cdot,0) = a$  in $L^2$-sense.  Moreover, for sufficiently smooth functions 
$a$ and $f$, the solution $u$ to \eqref{(2.3)} satisfies 
the initial condition in a usual sense: $\lim_{t\to 0} 
u(\cdot,t) = a$ in $L^2(\OOO)$ (see Lemma 4 in Section 3 of \cite{LYP1}).
Consequently, the third line of \eqref{(2.3)} can be  interpreted as 
a generalized initial condition.


\section{Well-posedness results} 
\label{sec5}

\setcounter{section}{3}
\setcounter{equation}{0}

In this section, we demonstrate well-posedness of the initial-boundary value 
problem \eqref{(2.3)} with a semilinear term $f$ that satisfies the conditions 
\eqref{(3.3)} formulated below. 

Before starting with formulation of our findings, 
we introduce some notations and known results needed for further discussions. 

For an arbitrarily chosen constant $c_0>0$,
we define the elliptic operator $A_0$ as follows:
\begin{equation}
\label{(3.2)}
\left\{ \begin{array}{rl}
& (-A_0v)(x) := \sumij \ppp_i(a_{ij}(x)\ppp_jv(x)) - c_0v(x), \quad
x\in \OOO, \\
& \mathcal{D}(A_0) = \left\{ v\in H^2(\OOO);\,
\NUNU v + \sigma v = 0 \quad \mbox{on } \ppp\OOO \right\}.
\end{array}\right.
\end{equation}

We recall that in \eqref{(3.2)}, $\sigma \in C^{\infty}(\ppp\OOO)$, 
$\sigma(x)\ge 0,\ x\in \ppp\OOO$, and 
the coefficients $a_{ij}$ satisfy the conditions 
\eqref{(1.2)}.

Henceforth, by $\Vert \cdot\Vert$ and $(\cdot,\cdot)$ we denote the 
standard norm and
the scalar product  in $L^2(\OOO)$, respectively. 
It is well-known that the operator $A_0$ is self-adjoint and its resolvent 
is a compact operator.  Moreover, for a sufficiently large constant
$c_0>0$, Lemma 6 from \cite{LYP1} ensures  that $A_0$ is 
positive definite.
Therefore, by choosing the  constant $c_0>0$ large enough,
the spectrum of $A_0$ consists entirely of discrete positive 
eigenvalues $0 < \la_1 \le \la_2 \le \cdots$,
which are numbered according to their multiplicities and  
$\la_n \to \infty$ as $n\to \infty$.
Let $\va_n$ be an eigenvector corresponding to the eigenvalue $\la_n$ 
such that $A\va_n = \la_n\va_n$ and 
$(\va_n, \va_m) = 0$ if $n \ne m$ and $(\va_n,\va_n) = 1$. 
Then the system $\{ \va_n\}_{n\in \N}$ of the 
eigenvectors forms an
orthonormal basis in $L^2(\OOO)$ and for any $\gamma\ge 0$ we can define
the fractional powers $A_0^{\gamma}$ of the operator $A_0$ 
by the following relation (see, e.g., \cite{Pa}):
$$
A_0^{\gamma}v = \sum_{n=1}^{\infty} \la_n^{\gamma} (v,\va_n)\va_n,
$$
where
$$
v \in \mathcal{D}(A_0^{\gamma})
:= \left\{ v\in L^2(\OOO): \thinspace
\sum_{n=1}^{\infty} \la_n^{2\gamma} (v,\va_n)^2 < \infty\right\}
$$
and 
$$
\Vert v\Vert_{\DDD(A_0^{\gamma})}
:= \Vert A_0^{\gamma}v\Vert = \left( \sum_{n=1}^{\infty}
\la_n^{2\gamma} (v,\va_n)^2 \right)^{\frac{1}{2}}.
$$
We also mention the important inclusion $\mathcal{D}(A_0^{\gamma}) \subset 
H^{2\gamma}(\OOO)$, where $H^{2\gamma}(\OOO)$ is the Sobolev-Slobodeckij space.

Now we introduce the conditions posed on the semilinear term 
$f$ of the equation from the problem \eqref{(2.3)}. 

For a constant $\gamma \in \left( \frac{3}{4},\, 1 \right)$ and a constant $m>0$, there exists a constant $C_f=C_f(m)>0$ such 
that
\begin{equation}
\label{(3.3)}
\left\{ \begin{array}{rl}
& \mbox{(i)} \Vert f(v)\Vert \le C_f, \quad 
\Vert f(v_1) - f(v_2)\Vert \le C_f\Vert v_1-v_2\Vert_{\DDDA} \\ 
& \mbox{if } \Vert v\Vert_{\DDD(A_0^{\gamma})}, 
\Vert v_1\Vert_{\DDD(A_0^{\gamma})}, 
\Vert v_2\Vert_{\DDD(A_0^{\gamma})}\le m \\
& \mbox{and}\\
& \mbox{(ii) there exists a constant\ } \ep > 0
\mbox{ such that }\\
& \Vert f(v)\Vert_{H^{2\ep}(\OOO)} \le C_f(m) \quad
\mbox{if } \Vert v\Vert_{\mathcal{D}(A_0^{\gamma})} \le m.
\end{array}\right.
\end{equation}
Henceforth, by $C>0$, $C_0, C_1 > 0$, etc., we denote generic constants, which 
are independent of the functions $u, v$, etc. under consideration, and 
we write 
$C_f$, $C(m)$ in the case where 
we need to specify a dependence on related quantities.
Before we state and prove the main results of this section, let us discuss 
some examples of semilinear functions, which satisfy the condition 
\eqref{(3.3)}. 

\begin{example}
\label{e1}
For $f\in C^1(\R)$, by setting $f(v):= f(v(x))$ for $x \in \OOO$,
we define 
$f: \DDD(A_0^{\gamma})\, \longrightarrow L^2(\OOO)$, 
$\frac{3}{4} < \gamma < 1$.
\end{example} 

Let us verify that the function from Example \ref{e1} satisfies the conditions 
\eqref{(3.3)}. For the spatial dimensions 
$d \le 3$, the Sobolev embedding leads to the inclusions
$\DDDA \subset H^{2\gamma}(\OOO) \subset C(\ooo{\OOO})$.  
Therefore, $\Vert v\Vert_{\DDD(A_0^{\gamma})}\le m$ yields
$\vert v(x)\vert \le C_0m$ for $x \in \OOO$.  Hence, by using the mean
value theorem, the first condition in \eqref{(3.3)} is satisfied if
we choose $C_f(m) = \Vert f\Vert_{C^1[-C_0m,\,C_0m]}$.

Next we have to verify the second condition in \eqref{(3.3)}.
To this end, we prove the following lemma:

\begin{lemma}
\label{l11.1}
Let $f \in C^1[-C_0m,C_0m]$ and $0 < \ep < \frac{3}{4}$.
Then
$$
\Vert f(w)\Vert_{\DDD(A_0^{\ep})} \le C(1+m)\Vert f\Vert_{C^1[-C_0m,C_0m]}
\quad \mbox{if } \Vert w\Vert_{\DDDA} \le m.
$$
\end{lemma}

\begin{proof}
Under the condition $0 < \ep < \frac{3}{4}$, the space
$\DDD(A_0^{\ep})$  coincides with the space $ H^{2\ep}(\OOO)$
(\cite{Fu}, \cite{Gri}).  
Then, for $f \in C^1[-C_0m,C_0m]$, the following norm estimates
in the Sobolev-Slobodeckij space $H^{2\ep}(\OOO)$ hold true:
\begin{align*}
& \Vert f(w)\Vert^2_{H^{2\ep}(\OOO)}
= \Vert f(w)\Vert_{L^2(\OOO)}^2
+ \int_{\OOO}\int_{\OOO} \frac{\vert f(w(x)) - f(w(y))\vert^2}
{\vert x-y\vert^{d+4\ep}} dxdy\\
\le& \Vert f\Vert^2_{C^1[-C_0m,C_0m]}\vert \OOO\vert^2
+ \int_{\OOO}\int_{\OOO} \frac{\Vert f\Vert^2_{C^1[-C_0m,C_0m]}
\vert w(x) - w(y)\vert^2}{\vert x-y\vert^{d+4\ep}} dxdy\\
\le &C(1 + \Vert w\Vert^2_{H^{2\ep}(\OOO)})\Vert f\Vert^2_{C^1[-C_0m,C_0m]}
\le C(1+m^2)\Vert f\Vert^2_{C^1[-C_0m,C_0m]}.
\end{align*}
\end{proof}
$\square$

The result formulated in Lemma \ref{l11.1} ensures the validity of the second  
condition from \eqref{(3.3)} for the function defined as in Example \ref{e1}.

\begin{example}
\label{e2}
We set
$$
f(v)(x):= \sum_{k=1}^d \mu_k(x)v(x)\ppp_kv(x),\quad x\in \OOO
$$
with some $\mu_k \in C^1(\ooo{\OOO})$ for $k=1,2,..., d$.  
In particular, a semilinear term of this type is contained in the 
time-fractional Burgers equation $\pppa u = \ppp_x^2u - 
u\ppp_xu$.
\end{example}

Let us verify that  the function defined in Example \ref{e2} satisfies the 
conditions \eqref{(3.3)}. 

For $\frac{3}{4} < \gamma < 1$ and a function 
$v \in \mathcal{D}(A_0^{\gamma})$,we immediately  get the norm estimates
$\Vert \nabla v\Vert \le C\Vert v\Vert_{H^1(\OOO)}
\le C\Vert A_0^{\gamma}v\Vert$.
Moreover, in the dimensions $d=1,2,3$,  the Sobolev embedding implies 
the inequality 
$\Vert v\Vert_{C(\ooo{\OOO})} \le C\Vert v\Vert_{\DDDA}$.

Applying the relations derived above, we obtain the inequalities from the 
first line in  \eqref{(3.3)}:
\begin{align*}
& \Vert f(v)\Vert \le C\sum_{k=1}^d \Vert v\ppp_kv\Vert 
\le C\Vert v\Vert_{C(\ooo{\OOO})}\sum_{k=1}^d \Vert \ppp_kv\Vert\\
\le& C\Vert v\Vert_{C(\ooo{\OOO})}\Vert v\Vert_{\DDDA}
\le C\Vert v\Vert^2_{\DDDA} \le Cm^2 =: C_f(m)
\end{align*}
and
\begin{align*}
& \Vert f(v_1)-f(v_2)\Vert 
= \left\Vert \sum_{k=1}^d \mu_k(v_1-v_2)\ppp_kv_1
+ \mu_kv_2\ppp_k(v_1-v_2)\right\Vert\\
\le& C( \Vert v_1-v_2\Vert_{C(\ooo{\OOO})}\Vert \nabla v_1\Vert
+ \Vert v_2\Vert_{C(\ooo{\OOO})}\Vert \nabla (v_1-v_2)\Vert)\\
\le &C\max\{ \Vert v_1\Vert_{\DDDA},\, \Vert v_2\Vert_{\DDDA}\}
\Vert v_1 - v_2\Vert_{\DDDA}.
\end{align*}

Now we verify the second condition from \eqref{(3.3)}.  Since $\mu_k\in 
C^1(\ooo{\OOO})$, it suffices to prove that we can choose some constant
$C(m) > 0$ such that  
\begin{equation}
\label{(11.1)}
\Vert v\ppp_kv\Vert_{H^{2\ep}(\OOO)} \le C(m), \quad 
k=1,..., d \quad \mbox{if } \Vert v\Vert_{\DDD(A_0^{\gamma})} \le m.
\end{equation}

For $\gamma > \frac{3}{4}$ and $\theta < 2\gamma - \frac{3}{2}$, the Sobolev 
embedding
\begin{equation}
\label{(11.2)}
H^{2\gamma}(\OOO) \subset C^{\theta}(\ooo{\OOO})   
\end{equation}
holds true in the dimensions $d=1,2,3$ (see, e.g., 
Theorem 1.4.4.1 (p. 27) in \cite{Gri}).
In \eqref{(11.2)}, $C^{\theta}(\ooo{\OOO})$ denotes the Schauder space of 
uniform H\"older continuous functions on $\ooo{\OOO}$, which is defined 
  as follows (see, e.g., \cite{GT},
\cite{LU}):  A function $w$ is said to belong to the space 
$C^{\theta}(\ooo{\OOO})$ if 
$$
\sup_{x, x'\in \OOO, \, x \ne x'} 
\frac{\vert w(x) - w(x')\vert}{\vert x-x'\vert^{\theta}}
< \infty.
$$  
For  $w \in C^{\theta}(\ooo{\OOO})$, we define the norm
$$
\Vert w\Vert_{C^{\theta}(\ooo{\OOO})}
:= \Vert w\Vert_{C(\ooo{\OOO})}
+ \sup_{x, x'\in \OOO, \, x \ne x'}  
\frac{\vert w(x) - w(x')\vert}{\vert x-x'\vert^{\theta}}.
$$

In what follows, we choose and fix a  number $\ep \in (0,1)$ that satisfies 
the condition
\begin{equation}
\label{(11.3)}
\ep < \min\left\{ \hhalf\theta, \, \frac{1}{4} \right\}.
\end{equation}

Now, using the inequalities
$$
\Vert v\Vert_{H^1(\OOO)} \le C\Vert v\Vert_{\DDD(A_0^{\gamma})}
\le Cm,
$$ 
we estimate the norm
\begin{align*}
& \Vert v\ppp_kv\Vert^2_{H^{2\ep}(\OOO)}
= \Vert v\ppp_k v\Vert^2_{L^2(\OOO)}
+ \int_{\OOO}\int_{\OOO}
\frac{\vert (v\ppp_k v)(x) - (v\ppp_k v)(y)\vert^2}
{\vert x-y\vert^{d+4\ep}} dxdy\\
=:& I_1+I_2.
\end{align*}
The inclusion \eqref{(11.2)} leads to the inequalities
\begin{align*}
&I_1 \le \Vert v\Vert^2_{C(\ooo{\OOO})}\Vert \ppp_kv\Vert^2_{L^2(\OOO)}
\le C\Vert v\Vert_{\DDD(A_0^{\gamma})}^2\Vert v\Vert_{H^1(\OOO)}^2\\
\le & C\Vert v\Vert^4_{\DDD(A_0^{\gamma})} \le Cm^4.
\end{align*}
As to the term $I_2$, we employ the estimates 
\begin{align*}
& \vert (v\ppp_kv)(x) - (v\ppp_kv)(y)\vert^2
\le (\vert v(x)(\ppp_kv(x) - \ppp_kv(y))\vert 
+ \vert v(x) - v(y)\vert \vert \ppp_kv(y)\vert)^2\\
\le & 2\Vert v\Vert^2_{C(\ooo{\OOO})}
\vert \ppp_kv(x) - \ppp_kv(y)\vert^2
+ 2\vert \ppp_kv(y)\vert^2 \vert v(x) - v(y)\vert^2
\end{align*}
and proceed as follows:
\begin{align*}
& I_2 \le C\Vert v\Vert_{C(\ooo{\OOO})}^2
\int_{\OOO}\int_{\OOO}
\frac{\vert (\ppp_k v)(x) - (\ppp_k v)(y)\vert^2}
{\vert x-y\vert^{d+4\ep}} dxdy\\
+ & C\int_{\OOO} \int_{\OOO} \vert \ppp_kv(y)\vert^2 
 \frac{\vert v(x) - v(y)\vert^2}{\vert x-y\vert^{d+4\ep}} dxdy\\
=: & I_{21} + I_{22}.
\end{align*}
As already mentioned, the Sobolev embedding yields the inequality 
$\Vert v\Vert_{C(\ooo{\OOO})}
\le C\Vert v\Vert_{\DDD(A_0^{\gamma})}$.
Moreover, since $0<\ep<\frac{1}{4}$ by the condition \eqref{(11.3)} and 
$\gamma > \frac{3}{4}$,  the relations $2\ep + 1 < 2\frac{1}{4} + 1 
= \frac{3}{2} < 2\gamma$ hold true. Thus, for these parameter values, 
we obtain the inequalities
$$
\Vert \ppp_kv\Vert_{H^{2\ep}(\OOO)} \le C\Vert v\Vert_{H^{2\ep+1}(\OOO)}
\le C\Vert v\Vert_{H^{2\gamma}(\OOO)} \le C\Vert v\Vert_{\DDD(A_0^{\gamma})}.
$$
Hence
$$
I_{21} \le C\Vert v\Vert^2_{C(\ooo{\OOO})}
\Vert \ppp_kv\Vert^2_{H^{2\ep}(\OOO)}
\le C\Vert v\Vert^2_{C(\ooo{\OOO})}
\Vert v\Vert^2_{\DDD(A_0^{\gamma})} \le Cm^4.
$$
To estimate the term $I_{22}$, we mention that the inclusion \eqref{(11.2)} 
yields the inequality
$$
\vert v(x) - v(y)\vert \le C\Vert v\Vert_{\DDD(A_0^{\gamma})}
\vert x-y\vert^{\theta}, \quad x,y \in \ooo{\OOO}.
$$
Thus we obtain the following estimate:
$$
I_{22} \le C\Vert v\Vert_{\DDD(A_0^{\gamma})}^2 \int_{\OOO}
\vert \ppp_kv(y)\vert^2\left(\int_{\OOO} \vert x-y\vert^{2\theta-d-4\ep}
dx \right) dy.
$$
Introducing the polar coordinate $r:= \vert x-y\vert$ for the inner integral 
(the one with the integration variable $x$) 
with a fixed $y\in \OOO$ and choosing a number $R_0>0$ sufficiently large such 
that 
$\max_{x,y\in \ooo{\OOO}} \vert x-y\vert < R_0$, we proceed as follows: 
\begin{align*}
& \int_{\OOO} \vert x-y\vert^{2\theta-d-4\ep}dx 
\le \int_{\vert x-y\vert < R_0} \vert x-y\vert^{2\theta-d-4\ep} dx\\
\le& C\int^{R_0}_0 r^{2\theta-d-4\ep}r^{d-1} dr
\le C\int^{R_0}_0 r^{(2\theta-4\ep)-1} dr =: C_1 < +\infty.
\end{align*}
To derive the last inequality, we used the condition $2\theta-4\ep > 0$ that 
holds true because of \eqref{(11.3)}.
Combining the results derived above, we get the desired estimate
$$
I_{22} \le CC_1\Vert v\Vert^2_{\DDD(A_0^{\gamma})}\int_{\OOO}
\vert \ppp_kv(y)\vert^2 dy
\le C\Vert v\Vert_{\DDD(A_0^{\gamma})}^2\Vert v\Vert^2_{H^1(\OOO)}
\le Cm^4.
$$
Thus, the conditions \eqref{(3.3)} are satisfied for the semilinear term 
$f$ defined as in Example \ref{e2}.

Now we proceed to some results regarding the local unique 
existence of a solution to the initial-boundary value problem \eqref{(2.3)}. 

\begin{theorem}
\label{t3.1}
Let a semilinear term $f$ satisfy the conditions \eqref{(3.3)} with $m>0$ and 
$\Vert a\Vert_{\DDDA} \le m$.  

Then there exists a constant $T=T(m)>0$
such that the initial-boundary value problem \eqref{(2.3)} 
possesses a unique solution $u=u_a(x,t)$ satisfying the inclusions 
\begin{equation}
\label{(3.4)}
u_a \in L^2(0,T;H^2(\OOO)) \cap C([0,T];\DDDA), \quad
u_a-a \in \HH(0,T;L^2(\OOO)).                      
\end{equation}

Moreover, there exists a constant $C(m)>0$, such that
\begin{equation}
\label{(3.5)}
\Vert u_a - u_b\Vert_{L^2(0,T;H^2(\OOO))}
\le C\Vert a-b\Vert_{\DDDA}                  
\end{equation}
provided that $\Vert a\Vert_{\DDDA}, \Vert b\Vert_{\DDDA} \le m$.
\end{theorem}

We recall that by $d=1,2,3$, the Sobolev embedding yields 
$\DDD(A_0^{\gamma}) \subset C(\ooo{\OOO})$ and (3.6) implies 
$u_a \in C(\ooo{\OOO}\times [0,T])$.

The results formulated in Theorem \ref{t3.1} are similar to the ones 
well-known for 
the partial differential equations of parabolic type that correspond to 
the case $\alpha=1$ in \eqref{(2.3)} (see, e.g.,  \cite{He} or \cite{Pa}).
One of the available and most used methodologies for the proof of 
Theorem \ref{t3.1} in the case $\alpha=1$ 
is the theory of the analytic semi-groups.  For the case $0<\alpha<1$, 
our proof relies on 
a similar idea, that is, on employing the solution operators $S(t)$ and 
$K(t)$ defined below 
even if the operator  $S(t)$ does not possess any semi-group property.
In several published works, 
a similar approach was employed for the equations with a self-adjoint operator 
$A$ in the case when the constant $C_f$ can be chosen independently of 
the constant $m>0$ (see, e.g., Section 6.4.1 of 
Chapter 6 in \cite{J}).  However, in general, 
the constant $C_f$ depends on $m$, in other words, the semilinear term 
$f$ is not globally Lipschitz continuous. Thus, the arguments employed 
in the earlier publications are not immediately applicable for our situation 
and in what follows we present their suitable modifications. 

\begin{proof}

First, we introduce the notations (\cite{GLY}, \cite{KRY})
\begin{equation}
\label{(5.1)}
S(t)a = \sum_{n=1}^{\infty} E_{\alpha,1}(-\la_n t^{\alpha})
(a,\va_n)\va_n, \quad a\in L^2(\OOO), \thinspace t>0,  
\end{equation}
\begin{equation}
\label{(5.2)}
K(t)a = -A_0^{-1}S'(t)a 
= \sum_{n=1}^{\infty} t^{\alpha-1}E_{\alpha,\alpha}(-\la_n t^{\alpha})
(a,\va_n)\va_n, \quad a\in L^2(\OOO), \thinspace t>0,
\end{equation}
where $E_{\alpha,\beta}(z)$ denotes the Mittag-Leffler function defined by 
a convergent series as follows:
$$
E_{\alpha,\beta}(z) = \sum_{k=0}^\infty \frac{z^k}{\Gamma(\alpha\, k + \beta)},
\ \alpha >0,\ \beta \in \mathbb{C},\ z \in \mathbb{C}.
$$
It follows directly from the definitions given above that
$A_0^{\gamma}K(t)a = K(t)A_0^{\gamma}a$
and $A_0^{\gamma}S(t)a = S(t)A_0^{\gamma}a$ for $a \in \mathcal{D}
(A_0^{\gamma})$.
Moreover, the inequality (see, e.g., Theorem 1.6 (p. 35) in \cite{Po})
$$
\max \{ \vert E_{\alpha,1}(-\la_nt^{\alpha})\vert, \, 
\vert E_{\alpha,\alpha}(-\la_nt^{\alpha})\vert \} 
\le \frac{C}{1+\la_nt^{\alpha}}
\quad \mbox{for all } t>0
$$
implicates the estimations (\cite{GLY})
\begin{equation}
\label{(5.3)}
\left\{ \begin{array}{l}
\Vert A_0^{\gamma}S(t)a\Vert \le Ct^{-\alpha\gamma}\Vert a\Vert, \\
\Vert A_0^{\gamma}K(t)a\Vert \le Ct^{\alpha(1-\gamma)-1}
\Vert a\Vert, \quad a \in L^2(\OOO), \thinspace t > 0, \thinspace
0 \le \gamma \le 1.
\end{array}\right.                    
\end{equation}
We also employ the operator
\begin{equation}
\label{(7.1)}
Q(t)u(t) := \sum_{j=1}^d b_j(\cdot,t)\ppp_ju(\cdot,t) 
+ (c(\cdot,t)+c_0)u(\cdot,t) \quad \mbox{in $\OOO$}, \, 0<t<T.   
\end{equation}
Then for any function $u\in \mathcal{D}(A_0)$, the norm estimates 
\begin{equation}
\label{(7.2)}
\Vert \Ahalf Q(t)u(t)\Vert \le C\Vert Q(t)u(t)\Vert_{H^1(\OOO)}
\le C\Vert Q(t)u(t)\Vert_{H^2(\OOO)} \le C\Vert A_0u(t)\Vert            
\end{equation}
hold true (\cite{GLY}, \cite{KRY}).

With the notations introduced above, the problem 
\eqref{(2.3)}  can be formally rewritten as an integral equation 
(\cite{LYP1})
\begin{equation}
\label{(7.3)}
u(t) = S(t)a + \int^t_0 K(t-s)Q(s)u(s) ds 
\end{equation}
$$
+ \int^t_0 K(t-s)f(u(s)) ds,
\quad 0<t<T.                                        
$$

By $\gamma \in \left( \frac{3}{4}, \, 1\right)$ 
in the condition \eqref{(3.3)} and a fixed initial value  
$a \in \DDD(A_0^{\gamma})$, 
we define an operator $L: L^2(0,T;L^2(\OOO))
\, \longrightarrow \,  L^2(0,T;L^2(\OOO))$ by 
$$
(Lu)(t) := S(t)a + \int^t_0 K(t-s)Q(s)u(s) ds + \int^t_0 K(t-s)f(u(s)) ds,
\quad 0<t<T.
$$
Choosing a constant $m>0$ arbitrarily, we set 
\begin{equation}
\label{(7.4)}
V:= \{ v\in C([0,T];\DDD(A_0^{\gamma}));\, 
\Vert u - S(\cdot)a\Vert_{C([0,T];\DDD(A_0^{\gamma}))} \le m\}.  
\end{equation}

Then we prove the following lemma:

\begin{lemma}
\label{l7.1}
Let $H \in C([0,T]; L^2(\OOO))$.  Then
$$
\int^t_0 A_0^{\gamma}K(t-s)H(s) ds \in C([0,T];L^2(\OOO)),
$$
that is,
$$
\int^t_0 K(t-s)H(s) ds \in C([0,T];\mathcal{D}(A_0^{\gamma})).
$$
\end{lemma}

\begin{proof}
Let $0< \eta < t \le T$.  We have the representation
\begin{align*}
& \int^t_0 A_0^{\gamma}K(t-s)H(s) ds - \int^{\eta}_0 A_0^{\gamma}K(\eta-s)
H(s) ds\\
=& \int^t_0 A_0^{\gamma}K(s)H(t-s) ds - \int^{\eta}_0 A_0^{\gamma}K(s)
H(\eta-s) ds\\
=& \int^t_{\eta} A_0^{\gamma}K(s)H(t-s) ds 
+ \int^{\eta}_0 A_0^{\gamma}K(s) (H(t-s) - H(\eta-s)) ds\\
=: &I_1 + I_2.
\end{align*}
For the first integral, by \eqref{(5.3)} and $\gamma < 1$, 
we have the relations
\begin{align*}
&\Vert I_1\Vert
\le C\int^t_{\eta} s^{\alpha(1-\gamma)-1} \max_{0\le s\le t}
\Vert H(t-s)\Vert ds\\
\le& C\Vert H\Vert_{C([0,T];\DDD(A_0^{\gamma}))}
\frac{t^{\alpha(1-\gamma)} - \eta^{\alpha(1-\gamma)}}
{\alpha(1-\gamma)} \, \longrightarrow \, 0
\end{align*}
as $\eta \uparrow t$.  

Next, we get the following norm estimates
\begin{align*}
& \Vert I_2\Vert
=\left\Vert \int^{\eta}_0 A_0^{\gamma}K(s) (H(t-s) - H(\eta-s)) ds\right\Vert\\
\le & C\int^{\eta}_0 s^{(1-\gamma)\alpha-1} 
\max_{0\le \eta \le t \le T} \Vert H(t-s) - H(\eta-s)\Vert ds.
\end{align*}
For $H\in C([0,T];L^2(\OOO))$, the function
$$
\vert s^{(1-\gamma)\alpha-1} \vert
\max_{0\le \eta \le t \le T} \Vert H(t-s) - H(\eta-s)\Vert 
$$
is an integrable function with respect to $s \in (0,\eta)$ and 
$$
\lim_{\eta \uparrow t} s^{(1-\gamma)\alpha-1} 
\max_{0\le \eta \le t \le T} \Vert H(t-s) - H(\eta-s)\Vert = 0
$$
for almost all $s\in (0,\eta)$.
Hence, the Lebesgue convergence theorem implies the relation 
$\lim_{\eta \uparrow t} \Vert I_2\Vert = 0$ and the proof of 
Lemma \ref{l7.1} is completed.
\end{proof}

Now we proceed to the proof of Theorem \ref{t3.1}.
In view of \eqref{(5.1)}, the inclusion $a \in \DDDA$ 
implies
\begin{equation}
\label{(7.5)}
S(t)a \in C([0,T];\DDDA).                       
\end{equation}
Indeed, 
\begin{align*}
& \Vert A_0^{\gamma}(S(t)a - S(s)a)\Vert^2
= \Vert S(t)(A_0^{\gamma}a) - S(s)(A_0^{\gamma}a)\Vert^2\\
=& \sum_{n=1}^{\infty} \vert E_{\alpha,1}(-\la_nt^{\alpha})
- E_{\alpha,1}(-\la_ns^{\alpha})\vert^2 \vert (A_0^{\gamma}a,\va_n)\vert^2.
\end{align*}
Applying the Lebesgue convergence theorem and the estimate (see, e.g., 
Theorem 1.6 (p. 35) in \cite{Po})
$$
\vert E_{\alpha,1}(-\la_nt^{\alpha})\vert \le \frac{C}{1+\la_nt^{\alpha}}
\quad \mbox{for all } n\in \N \mbox{ and } t>0
$$
we can verify the inclusion \eqref{(7.5)}.

Because of the condition \eqref{(3.3)} and  $\mathcal{D}(A_0^{\gamma})
\subset H^1(\OOO)$, 
for $v \in C([0,T];\DDD(A_0^{\gamma}))$, we get the inclusions  
$f(v) \in C([0,T];L^2(\OOO))$ and $Qv \in C([0,T];L^2(\OOO))$.
Now, applying Lemma \ref{l7.1}, in view of \eqref{(7.5)}, we arrive at the  
inclusion
\begin{equation}
\label{(7.6)}
Lv \in C([0,T]; \DDDA) \quad \mbox{for } v \in C([0,T]; \DDDA).
\end{equation}

For the further proof, we need
the following properties that
are valid for a sufficiently small $T>0$:
\\
(i) $LV \subset V$, $V$ being the set defined by \eqref{(7.4)}.
\\
(ii) There exists a constant $\rho \in (0,1)$ such that for any 
$u_1, u_2 \in V$, the norm estimate
$$
\Vert Lu_1 - Lu_2\Vert_{C[0,T];\DDDA)} 
\le \rho\Vert u_1 - u_2\Vert_{C[0,T];\DDDA)}, \quad 0<t<T
$$
holds true. 
\\
{\bf Proof of (i).}
\\
Let $u \in V$.  Then, the inclusion \eqref{(7.6)} implicates 
$Lu \in C([0,T];\DDDA)$.
Now we consider the expression 
\begin{equation}
\label{(7.7)}
A_0^{\gamma}(Lu(t) - S(t)a)
= \int^t_0 A_0^{\gamma}K(t-s)Q(s)u(s) ds + \int^t_0 A_0^{\gamma}K(t-s)
f(u(s)) ds.                                     
\end{equation}
For any $u \in V$, using the norm estimates 
$$
\Vert a\Vert_{\mathcal{D}(A_0^{\gamma})}
= \Vert A_0^{\gamma}a\Vert \le m, \quad
\Vert u-S(\cdot)a\Vert_{C([0,T];\mathcal{D}(A_0^{\gamma}))} \le m,
$$
we obtain
\begin{equation}
\label{(7.8)}
\Vert u(t)\Vert_{\DDDA} \le m + \Vert A_0^{\gamma}S(t)a\Vert 
= m + \Vert S(t)A_0^{\gamma}a\Vert 
\le m + C_1m =: C_2m.
\end{equation}
The first condition from \eqref{(3.3)} implies that
\begin{equation}
\label{(7.9)}
\Vert f(u(t)) \Vert \le C_f(C_2m) \quad \mbox{for all $u\in V$
and $0<t<T$.}    
\end{equation}
Furthermore, for $\gamma \ge \hhalf$, \eqref{(7.8)} leads 
to the norm estimates
\begin{equation}
\label{(7.10)}
\Vert Q(s)u(s)\Vert \le C\Vert u(s)\Vert_{H^1(\OOO)}
\le C_3\Vert A_0^{\hhalf}u(s)\Vert \le C_4\Vert u(t)\Vert_{\DDD(A_0^{\gamma})}
\le C_4C_2m.                         
\end{equation}
Applying the inequalities \eqref{(7.9)} and \eqref{(7.10)} for the expression 
\eqref{(7.7)}, by means of \eqref{(5.3)},
we obtain the norm estimates
\begin{align*}
& \Vert Lu(t) - S(t)a \Vert_{\DDDA} \\
= & \left\Vert \int^t_0 A_0^{\gamma}K(t-s)Q(s)u(s) ds
+ \int^t_0 A_0^{\gamma}K(t-s)f(u(s)) ds\right\Vert\\
\le & C\int^t_0 (t-s)^{(1-\gamma)\alpha-1} (C_2C_4m + C_f(C_2m))ds
\le C_5\frac{t^{(1-\gamma)\alpha}}{(1-\gamma)\alpha}
\le C_5\frac{T^{(1-\gamma)\alpha}}{(1-\gamma)\alpha}.
\end{align*}
The constant $C_5>0$ depends on $m>0$ but is independent 
on $T>0$.
Therefore, choosing $T>0$ sufficiently small, we complete the proof of 
the property (i).  $\square$
\\
{\bf Proof of (ii).}
\\
Estimate \eqref{(7.8)} yields that $\Vert u_1(t)\Vert_{\DDDA} \le C_2m$ and 
$\Vert u_2(t)\Vert_{\DDDA} \le C_2m$ for any $u_1, u_2 \in V$.
The condition \eqref{(3.3)} leads then to the norm estimate
$$
\Vert f(u_1(s)) - f(u_2(s))\Vert \le C_f(C_2m)\Vert u_1(s)-u_2(s)\Vert
_{\DDDA}, \quad 0<s<T.
$$
Hence, by \eqref{(7.10)}, we have the following chain of estimates:
\begin{align*}
& \Vert Lu_1(t) - Lu_2(t)\Vert_{\DDDA} \\
= & \biggl\Vert \int^t_0 A_0^{\gamma}K(t-s)Q(s)(u_1(s)-u_2(s)) ds \\
+  & \int^t_0 A_0^{\gamma}K(t-s)(f(u_1(s)) - f(u_2(s))) ds\biggr\Vert\\
\le & C\int^t_0 (t-s)^{\alpha(1-\gamma)-1}\Vert (u_1-u_2)(s)\Vert
_{\DDDA} ds\\
+& C_f(C_2m)\int^t_0 (t-s)^{\alpha(1-\gamma)-1}\Vert (u_1-u_2)(s)\Vert
_{\DDDA} ds\\
\le & C_6T^{\alpha(1-\gamma)}
\sup_{0<s<T} \Vert u_1(s) - u_2(s)\Vert_{\DDDA}.
\end{align*}
In the last inequality, the constant $C_6$ is independent of $T$, and thus 
we can  choose a sufficiently small constant $T>0$ satisfying 
the inequality
$$
\rho:= C_6T^{\alpha(1-\gamma)} < 1.
$$
The proof of the property (ii) is completed.
$\square$

Due to the properties (i) and (ii), the contraction theorem can be applied to 
the equation $u=Lu$.  As a result, this equation has a unique 
solution $u\in V$  
for $0<t<T$.  This solution $u\in C([0,T];\mathcal{D}(A_0^{\gamma}))$ 
satisfies the estimates 
\eqref{(7.8)} and the equation
\begin{equation}
\label{(7.11)}
u(t) = S(t)a + \int^t_0 K(t-s)Q(s)u(s) ds 
\end{equation}
$$
+ \int^t_0 K(t-s)f(u(s)) ds, \quad 0<t<T.          
$$

Next, for the solution $u$ of the equation $u=Lu$,  we prove the inclusions 
\eqref{(3.4)}.

In the condition \eqref{(3.3)}, we can choose a sufficiently small $\ep>0$ 
such that $0<\ep<\hhalf$.
By the equation \eqref{(7.11)}, we obtain
\begin{align*}
& A_0u(t) = A_0^{1-\gamma}S(t)A_0^{\gamma}a
+ \int^t_0 A_0^{1-\ep}K(t-s)A_0^{\ep}Q(s)u(s) ds\\
+& \int^t_0 A_0^{1-\ep}K(t-s)A_0^{\ep}f(u(s)) ds, \quad 0<t<T.
\end{align*}
In the last equation, $a\in \DDD(A_0^{\gamma})$ and the estimates 
$$
\Vert A_0^{\ep}Q(s)u(s)\Vert \le C\Vert A_0^{\hhalf}(Q(s)u(s))\Vert 
\le C\Vert u(s)\Vert_{H^2(\OOO)} \le C\Vert A_0u(s)\Vert
$$
hold true because of \eqref{(7.2)}.  Furthermore, by the second condition 
in \eqref{(3.3)}, the inequality \eqref{(7.8)} 
yields the estimate $\Vert A_0^{\ep}f(u(s))\Vert 
\le C_f(C_2m)$.  
Thus, we get the chain of the inequalities
\begin{align*}
& \Vert A_0u(t)\Vert
\le Ct^{-\alpha(1-\gamma)}\Vert A_0^{\gamma}a\Vert 
+ C\int^t_0 (t-s)^{\alpha\ep-1}\Vert A_0u(s)\Vert ds\\
+ & C\int^t_0 (t-s)^{\alpha\ep-1}C_f(C_2m) ds\\
\le & Ct^{-\alpha(1-\gamma)}\Vert A_0^{\gamma}a\Vert + C_f(C_2m)
+ C\int^t_0 (t-s)^{\alpha\ep-1}\Vert A_0u(s)\Vert ds, \quad
0<t<T.
\end{align*}
For $0<\alpha<1$, the inequality $-\alpha(1-\gamma) > -1$ and 
thus the inclusion
$t^{-\alpha(1-\gamma)} \in L^1(0,T)$ hold true. 
Application of the generalized Gronwall inequality yields the estimates
\begin{align*}
& \Vert A_0u(t)\Vert
\le (Ct^{-\alpha(1-\gamma)}\Vert A_0^{\gamma}a\Vert + C_f(C_2m)) \\ 
+ & C\int^t_0 (t-s)^{\alpha\ep-1}(s^{-\alpha(1-\gamma)}
\Vert A_0^{\gamma}a\Vert + C_f(C_2m)) ds\\
\le&  Ct^{-\alpha(1-\gamma)}\Vert A_0^{\gamma}a\Vert 
+ C_f(C_2m)
+ (\Vert A_0^{\gamma}a\Vert + C_f(C_2m))t^{\alpha(\ep-(1-\gamma))},
\quad 0<t<T.
\end{align*}
Therefore, noting that $-\alpha(1-\gamma) < \alpha(\ep-(1-\gamma))$, we have 
the norm estimate 
$$
\Vert A_0u(t)\Vert \le C_7(1+T^{\alpha\ep})(t^{-\alpha(1-\gamma)} + 1),
\quad 0<t<T,
$$
where $C_7>0$ depends on $\Vert A_0^{\gamma}a\Vert$ and the constants
$C_f$, $C_2$, $m$, $\alpha$, and $\ep$.
For $\hhalf < \gamma \le 1$, we can directly verify that 
$-2\alpha(1-\gamma) > -1$, 
so that 
$\int^T_0 \Vert A_0u(t)\Vert^2 dt < \infty$, that is, the inclusion
\begin{equation}
\label{(7.12)}
u \in L^2(0,T;H^2(\OOO))                         
\end{equation}
holds true. 

It remains to prove that $u-a \in H_{\alpha}(0,T;L^2(\OOO))$.  
The inequality \eqref{(7.9)} implies the inclusion 
$f(u) \in L^2(0,T;L^2(\OOO))$.

The estimate \eqref{(7.10)} and the inclusion \eqref{(7.12)} result in the 
inclusion $Qu \in L^2(0,T;L^2(\OOO))$.  
Then we apply  the estimate (ii) of Lemma 1 from \cite{LYP1} and obtain the 
inclusion
$$
\int^t_0 K(t-s)Q(s)u(s) ds + \int^t_0 K(t-s)f(u(s)) ds
\in \HH(0,T;L^2(\OOO)).
$$
Applying the estimate (iii) of Lemma 1 from \cite{LYP1} to 
the equation \eqref{(7.3)}, we arrive at 
the inclusion $u-a 
\in \HH(0,T;L^2(\OOO))$, which completes the proof of the relation 
\eqref{(3.4)} from the theorem.

Finally, we have to prove the estimate \eqref{(3.5)}. By the construction of 
the solutions $u_a, u_b$ as the fixed points of the equation $u=Lu$, we have 
the inequalities
\begin{equation}\label{(5.20)}
\Vert u_a\Vert_{L^2(0,T;H^2(\OOO))}  \le C(m), \,
\Vert u_b\Vert_{L^2(0,T;H^2(\OOO))} \le C(m).
\end{equation}
On the other hand, 
\begin{align*}
&u_a(t) - u_b(t)
= S(t)a - S(t)b + \int^t_0 K(t-s)Q(s)(u_a-u_b)(s) ds\\
+& \int^t_0 K(t-s)(f(u_a(s)) - f(u_b(s))) ds, \quad 0<t<T. 
\end{align*}
In view of \eqref{(5.20)}, we can use the condition \eqref{(3.3)} and apply
the generalized Gronwall inequality.
Further details of the derivations are similar to the ones employed
in the proof of Theorem 1 from \cite{LYP1} 
and we omit them here. 
Thus, the proof of Theorem \ref{t3.1} is completed.
\end{proof}


\section{Comparison principles} 
\label{sec6}

\setcounter{section}{4}
\setcounter{equation}{0}

In this section, we derive some comparison principles for solutions to the 
initial-boundary value problems \eqref{(2.3)} for the semilinear
time-fractional diffusion equations.  

In that follows, we suppose that  
the semilinear term $f$ from \eqref{(2.3)}  depends only on the spatial 
variable $x$ and the unknown function $u$,
but not on gradient of $u$. Moreover, we introduce a class of semilinear 
terms $f$  via the smooth functions from the space
$C^1(\ooo{\OOO}\times [-m,\, m])$.  
For a function $\mathfrak{f} \in C^1(\ooo{\OOO}\times [-m,\, m])$, we define 
a mapping
$f: \{v\in \DDD(A_0^{\gamma}):\, \Vert v\Vert_{\DDD(A_0^{\gamma})} \le m\} 
\, \longrightarrow \, L^2(\OOO)$ by 
\begin{equation}
\label{(3.6a)}
f(v):= \mathfrak{f}(x,\, v(x)), \quad x\in \OOO,\, 0<t<T.      
\end{equation}
In the rest of this section, the semilinear term $f$ is determined by the 
function $\mathfrak{f}$ through the relation 
\eqref{(3.6a)}.
For a fixed constant $M>0$, we set 
\begin{equation}
\label{(3.7a)}
\mathcal{F}_M:= \{ f\in C^1(\ooo{\OOO}\times [-m,\,m]);\,
\Vert f\Vert_{C^1(\ooo{\OOO}\times [-m,\,m])} \le M\}.      
\end{equation}

Now we  formulate and prove the first comparison principle for 
the initial-boundary value problems  \eqref{(2.3)} for the semilinear 
time-fractional 
diffusion equations.

\begin{theorem}
\label{t3.2}
For $f_1, f_2 \in \mathcal{F}_M$ and $a_1, a_2 \in \DDD(A_0^{\gamma})$,  
we assume that there exist solutions $u(f_k,a_k)$, $k=1,2$ to the 
initial-boundary value problem \eqref{(2.3)}
with 
the semilinear terms $f_k$, $k=1,2$ and the initial values $a_k$, $k=1,2$, 
respectively, which 
satisfy the inclusions \eqref{(3.4)} and the inequalities
\begin{equation}
\label{(3.8)}
\vert u(f_k,a_k)(x,t)\vert \le m, \quad x\in \OOO,\, 
0<t<T, \, k=1,2.                              
\end{equation}
If $f_1(\cdot,\cdot) \ge f_2(\cdot,\cdot)$ on $\OOO \times (-m, m)$ and 
$a_1(\cdot)\ge a_2(\cdot)$ in $\OOO$, then 
\begin{equation}
\label{(3.9)}
u(f_1,a_1)(x,t) \ge u(f_2,a_2)(x,t) \quad \mbox{in } \OOO \times (0,T).
\end{equation}
\end{theorem}

\begin{proof}

Henceforth, for simplicity, we denote $f_k(x,u_j(x,t))$ by 
$f_k(u_j)$ and
$u(f_k,a_k)$ by $u_k$ for $j,k=1,2$.  

For $d=1,2,3$ and $\gamma > \frac{3}{4}$, application of the Sobolev embedding 
yields $\DDD(A_0^{\gamma}) \subset H^{2\gamma}(\OOO) \subset C(\ooo{\OOO})$.
Hence, \eqref{(3.4)} implies
\begin{equation}
\label{(8.1)}
u_1, u_2 \in C(\ooo{\OOO} \times [0,T]).     
\end{equation}
On the other hand, we have the representation
\begin{align*}
& f_1(u_1(x,t)) - f_2(u_2(x,t))
= f_1(u_1(x,t)) - f_1(u_2(x,t)) + (f_1-f_2)(u_2(x,t))\\
=& g(x,t)(u_1(x,t) - u_2(x,t)) + H(x,t), \quad
(x,t) \in \QQQQ,
\end{align*}
where we set
$$
g(x,t) := 
\left\{ \begin{array}{rl}
& \frac{f_1(u_1(x,t)) - f_1(u_2(x,t))}
{u_1(x,t) - u_2(x,t)} \quad \mbox{if } u_1(x,t) \ne u_2(x,t),\\
& f_1'(u_1(x,t)) \quad \mbox{if } u_1(x,t) = u_2(x,t),
\end{array}\right.
$$
and
$$
H(x,t) := (f_1-f_2)(u_2(x,t)), \quad (x,t) \in \QQQQ.
$$
By the inclusions \eqref{(8.1)} and $f \in C^1(\ooo{\OOO} \times [-m,m])$, 
we can verify that $g, H \in C(\ooo{\OOO}\times [0,T])$.  
Furthermore $f_1\ge f_2$ implies $H \ge 0$ in $\QQQQ$.

Setting now $y:= u_1 - u_2$ and $a:= a_1-a_2$, the function $y$ is a solution 
to the following initial-boundary value problem:
\begin{equation}
\label{(8.2)}
\left\{ \begin{array}{rl}
& \pppa (y-a) + Ay - g(x,t)y = H \ge 0 \quad 
\mbox{in } \QQQQ, \\
& \NUNU y + \sigma y = 0 \quad \mbox{on } \ppp\OOO \times (0,T), \\
& y(x,\cdot) - a(x) \in \HH(0,T) \quad \mbox{for almost all } x\in \OOO.
\end{array}\right.
\end{equation}
Since $g \in C(\ooo{\OOO}\times [0,T])$ does not in general satisfy 
the regularity condition \eqref{(1.2)}, we cannot directly apply Theorem 
2 from \cite{LYP1}.
To go forward,  we first approximate the function $g$ by 
$g_n \in C^{\infty}(\ooo{\OOO}\times
[0,T])$, $n\in \N$ such that $g_n \longrightarrow g$ in 
$C(\ooo{\OOO}\times [0,T])$ as $n\to \infty$ and consider a sequence of 
the following problems:
\begin{equation}
\label{(8.3)}
\left\{ \begin{array}{rl}
& \pppa (y_n-a) + Ay_n - g_n(x,t)y_n = H \quad 
\mbox{in } \QQQQ, \\
& \NUNU y_n + \sigma y_n = 0 \quad \mbox{on } \ppp\OOO \times (0,T), \\
& y_n(x,\cdot) - a(x) \in \HH(0,T) \quad \mbox{for almost all } x\in \OOO.
\end{array}\right.
\end{equation}
Then Theorem 1 from \cite{LYP1} yields that for any $n\in \N$ there exists a 
unique solution 
$y_n \in L^2(0,T;H^2(\OOO))$ to the problem \eqref{(8.3)} such that 
$y_n-a \in \HH(0,T;L^2(\OOO))$.  Moreover, since 
$g_n$ satisfies the regularity condition \eqref{(1.2)}, Theorem 2 from 
\cite{LYP1}
ensures that 
\begin{equation}
\label{(8.4)}
y_n(x,t)\ge 0 \quad \mbox{in } \QQQQ \mbox{ for each } n\in \N.
\end{equation}
On the other hand, setting $z_n:= y_n-y$ in $\QQQQ$, the equations 
\eqref{(8.2)} and 
\eqref{(8.3)} allow to characterize $z_n$ as solution to the problem
$$
\left\{ \begin{array}{rl}
& \pppa z_n + Az_n = (g_n-g)y + g_nz_n \quad 
\mbox{in } \QQQQ, \\
& \NUNU z_n + \sigma z_n = 0 \quad \mbox{on } \ppp\OOO \times (0,T), \\
& z_n(x,\cdot)  \in \HH(0,T) \quad \mbox{for almost all } x\in \OOO.
\end{array}\right.
$$

Similarly to the procedure that was employed at the beginning of the proof of 
Theorem 1 from \cite{LYP1}, we rewrite the first 
equation from the problem formulated above in the form
\begin{align*}
&\pppa z_n + A_0z_n\\
= &(g_n-g)y + g_nz_n + \sum_{j=1}^d b_j(t)\ppp_jz_n
+ (c_0+c(t))z_n \quad \mbox{in $\QQQQ$}
\end{align*}
and  obtain 
the  following chain of the estimates:
\begin{align*}
& \Vert A_0^{\hhalf} z_n(t)\Vert
= \biggl\Vert \int^t_0 A_0^{\hhalf}K(t-s)(g_n-g)(s)y(s) ds\\
+& \int^t_0 A_0^{\hhalf}K(t-s)
\left( (g_n(s)+c_0+c(s))z_n(s) 
+ \sum_{j=1}^d b_j(s)\ppp_jz_n(s) \right) ds \biggr\Vert\\
\le& C\int^t_0 (t-s)^{\hhalf\alpha-1} \Vert g_n-g\Vert
_{C(\ooo{\OOO}\times [0,T])}\Vert y(s)\Vert ds\\
+& \int^t_0 (t-s)^{\hhalf\alpha-1} (\Vert g_n\Vert_{C(\ooo{\OOO}\times [0,T])}
+ c_0 + \Vert c\Vert_{C(\ooo{\OOO}\times [0,T])} \\
+ & \max_{1\le j \le d} \Vert b_j\Vert_{C(\ooo{\OOO}\times [0,T])})
\Vert z_n(s)\Vert_{H^1(\OOO)} ds \\
\le & Ct^{\hhalf\alpha} \Vert g_n-g\Vert_{C(\ooo{\OOO}\times [0,T])}
+ C\int^t_0 (t-s)^{\hhalf\alpha-1} \Vert A_0^{\hhalf}z_n(s)\Vert ds\\
\le & C_1\Vert g_n-g\Vert_{C(\ooo{\OOO}\times [0,T])}
+ C_1\int^t_0 (t-s)^{\hhalf\alpha-1} \Vert A_0^{\hhalf}z_n(s)\Vert ds,
\quad 0<t<T.
\end{align*}
where the constant $C_1>0$ depends  on the norm 
$\Vert y\Vert_{C([0,T];L^2(\OOO))}$.
In the derivations presented above, we used the relation 
$\sup_{n\in \N}\Vert g_n\Vert_{C(\ooo{\OOO}\times [0,T])} < \infty$,
which is justified by  
$g_n \longrightarrow g$ in $C(\ooo{\OOO}\times [0,T])$ and 
$y:= u_1-u_2 \in C(\ooo{\OOO}\times [0,T])$ that follows from 
\eqref{(8.1)}.

Therefore, the generalized Gronwall inequality yields
\begin{align*}
& \Vert A_0^{\hhalf} z_n(t)\Vert
\le C_2\Vert g_n - g\Vert_{C(\ooo{\OOO}\times [0,T])}
+ C_2\int^t_0 (t-s)^{\hhalf\alpha-1} 
\Vert g_n - g\Vert_{C(\ooo{\OOO}\times [0,T])} ds\\
\le &C_3\Vert g_n - g\Vert_{C(\ooo{\OOO}\times [0,T])},
\quad 0<t<T,\, n\in \N.
\end{align*}
Hence, $z_n \longrightarrow 0$ in $L^{\infty}(0,T;\DDD(A_0^{\hhalf}))$
as $n\to \infty$, and thus  
$y_n \longrightarrow y$ in $L^{\infty}(0,T;\DDD(A_0^{\hhalf}))$
as $n\to \infty$.  Because of the relation
\eqref{(8.4)}, we readily get the inequality $y(x,t)\ge 0$ in $\QQQQ$.
Thus, the proof of Theorem \ref{t3.2} is completed.
\end{proof}

It is worth mentioning that in the formulation of Theorem \ref{t3.2}, 
the boundedness condition 
\eqref{(3.8)} has been assumed. However, this condition is not automatically 
guaranteed for solutions to the initial-boundary value problem \eqref{(2.3)}. 
In the rest of this section, we present another comparison principle in terms 
of the upper and lower solutions, which does not require this 
boundedness condition. 

\begin{definition}
\label{duls}
Let $f \in \mathcal{F}_M$.
The functions $\UPP$ and $\LOP$ satisfying \eqref{(3.4)} are 
called an upper solution and a lower solution for the solution 
$u$ to the problem \eqref{(2.3)}, respectively, 
if
\begin{equation}
\label{(3.10)}
\left\{ \begin{array}{rl}
& \pppa (\UPP - \ooo{a}) + A\UPP \ge f(\UPP) \quad 
\mbox{in } \OOO\times (0,T),               \\
& \ooo{a}(x) \ge u(x,0), \quad x\in \OOO,\\
& \NUNU \UPP + \sigma\UPP = 0 \quad \mbox{on } \ppp\OOO\times (0,T),\\
& \UPP - \ooo{a} \in \HH(0,T;L^2(\OOO)),
\end{array}\right.
\end{equation}
and
\begin{equation}
\label{(3.11)}
\left\{ \begin{array}{rl}
& \pppa (\LOP - \underline{a}) + A\LOP \le f(\LOP) \quad \mbox{in } \OOO\times 
(0,T),\\
& \underline{a}(x) \le u(x,0), \quad x\in \OOO,\\
& \NUNU \LOP + \sigma\LOP = 0 \quad \mbox{on } \ppp\OOO\times (0,T),\\
& \LOP - \underline{a} \in \HH(0,T;L^2(\OOO)),
\end{array}\right.
\end{equation}
with $\ooo{a}, \underline{a} \in L^2(\OOO)$.
\end{definition}

We recall that $\mathcal{F}_M$ is defined by \eqref{(3.7a)}.

Then the following result holds true:
\begin{theorem}
\label{t3.3}
For an arbitrarily chosen $T>0$,  we assume that there exist 
an upper solution $\ooo{u}$ and a lower solution $\underline{u}$
 to the problem \eqref{(2.3)} and that  
$\vert \UPP(x,t)\vert \le m$, $\vert \LOP(x,t)\vert \le m$ for $x\in \OOO$ and 
$0<t<T$.

Then there exists a unique solution $u_a$ to the problem \eqref{(2.3)} in 
the class \eqref{(3.4)} and
$$
\LOP(x,t) \le u_a(x,t) \le \UPP(x,t), \quad x\in \OOO,\, 0<t<T.
$$
\end{theorem}

\begin{proof}

For the functions $\UPPY, \LOWY \in L^2(0,T;L^2(\OOO))$ satisfying 
the condition 
$\LOWY(x,t) \le \UPPY(x,t)$ for $(x,t) \in \QQQQ$, we set 
$$
[\LOWY,\, \UPPY]:= \{ y\in L^2(0,T;L^2(\OOO));\, 
\LOWY(x,t) \le y(x,t) \le \UPPY(x,t)\},
$$
where the inequalities hold true for almost all  $(x,t) \in \QQQQ$. 

The key element of our proof of Theorem \ref{t3.3} is the following fixed 
point theorem in an ordered Banach space (see, e.g., \cite{Am}):

\begin{lemma}
\label{l9.1}
Let an operator $L: [\LOWY,\, \UPPY] \, \longrightarrow [\LOWY,\, \UPPY]
\subset L^2(0,T;L^2(\OOO))$ be compact and increasing, that is, 
$Lv \ge Lw$ in $\QQQQ$ if $v, w\in [\LOWY,\, \UPPY]$ and $v\ge w$ in 
$\QQQQ$.

Then the operator $L$ possesses fixed points $y^*$, $y_* \in 
[\LOWY,\, \UPPY]$ such that
$$
y_* = \lim_{k\to\infty} L^k\LOWY, \quad
y^* = \lim_{k\to\infty} L^k\UPPY \quad \mbox{in } L^2(0,T;L^2(\OOO)).
$$
\end{lemma}

Lemma \ref{l9.1} readily implies the following result:

\begin{lemma}
\label{l9.2}
In Lemma \ref{l9.1}, we further assume that the operator $L$ possesses 
a unique fixed point.
Then
$$
\LOWY \le y \le \UPPY \quad \mbox{in } \QQQQ.
$$
\end{lemma}

{\bf Proof of Lemma \ref{l9.2}.}
Indeed, since $L[\LOWY,\, \UPPY] \subset [\LOWY,\,\UPPY]$, we have the inequalities 
$\LOWY \le L\LOWY$ and $L\UPPY \le \UPPY$.
Moreover, since the operator $L$ is increasing, the inequality $L\LOWY \le
L\UPPY$ holds true.  Hence,
$$
\LOWY \le L\LOWY \le L\UPPY \le \UPPY \quad \mbox{in } \QQQQ.
$$
Applying the operator $L$ once again, we obtain
$$
\LOWY \le L\LOWY \le L^2\LOWY \le L^2\UPPY \le L\UPPY \le \UPPY.
$$
Repeating the same arguments, we get the inequalities 
$$
\LOWY \le L^k\LOWY \le L^k\UPPY \le \UPPY \quad \mbox{in } \QQQQ \mbox{ for }
k\in \N.
$$
By the uniqueness of the fixed point, letting $k \to \infty$, we finally get $\LOWY \le y \le \UPPY$.  
\\
$\square$
\\

Now we proceed to the proof of Theorem \ref{t3.3} and 
 define the operator $L$ from $[\UPP,\, \LOP] \subset  
L^2(0,T;L^2(\OOO))$ to $L^2(0,T;L^2(\OOO))$  as follows: 
For $u\in [\underline{u}, \, \ooo{u}]$, the function $v:=Lu$ is a solution to the problem
\begin{equation}
\label{(9.1)}
\left\{ \begin{array}{rl}
& \pppa (v-a) + Av + (M+1)v = (M+1)u + f(u) \quad \mbox{in } \QQQQ,\\
& v - a \in \HH(0,T;L^2(\OOO)), \\
& \NUNU v + \sigma v = 0 \quad \mbox{on } \ppp\OOO \times (0,T),
\end{array}\right.
\end{equation}
where $M>0$ is the constant in the definition \eqref{(3.7a)} of the set 
$\mathcal{F}_M$. 

First we verify that the operator $L$ 
is well-defined.
Indeed, let $u \in [\underline{u},\, \ooo{u}]$.
Because  $\vert \underline{u}(x,t)\vert \le m$ and $\vert \ooo{u}(x,t)\vert 
\le m$ for all $x\in \ooo{\OOO}$ and $0\le t \le T$, we obtain the estimate
$\vert u(x,t)\vert \le m$ for almost all $x\in \OOO$ and 
$t\in (0,T)$.  Therefore, the semilinear term $f(u(x,t))$ is defined for 
almost all $(x,t) \in \OOO \times (0,T)$ and $f(u(x,t)) 
\in L^{\infty}(\OOO\times (0,T))$.  Since $(M+1)u + f(u) 
\in L^{\infty}(\OOO\times (0,T)) \subset L^2(0,T;L^2(\OOO))$, we can apply 
Theorem 1 from \cite{LYP1} to the problem \eqref{(9.1)}. Thus, there exists 
a unique solution 
$v \in L^2(0,T;H^2(\OOO))$ to the problem \eqref{(9.1)} satisfying 
the inclusion
$v-a \in H_{\alpha}(0,T;L^2(\OOO))$. Moreover, the inequalities 
\begin{align*}
& \Vert (M+1)u+f(u)\Vert_{L^{\infty}(\OOO\times (0,T))} \\
\le & (M+1)\Vert u\Vert_{L^{\infty}(\OOO\times (0,T))}
+ \Vert f(u)\Vert_{L^{\infty}(\OOO\times (0,T))}
\le m(M+1)+M
\end{align*} 
lead to the norm estimates
\begin{equation}\label{(6.12)}
\Vert v-a\Vert_{H_{\alpha}(0,T;L^2(\OOO))} 
+ \Vert v\Vert_{L^2(0,T;H^2(\OOO))}
\end{equation}
\begin{align*}
\le& C(\Vert a\Vert_{H^1(\OOO)} + \Vert (M+1)u+f(u)\Vert
_{L^2(0,T;L^2(\OOO))}\\
\le& C(\Vert a\Vert_{H^1(\OOO)} 
+ ((M+1)m+M)T^{\hhalf}\vert \OOO\vert^{\hhalf})
=: C_4.
\end{align*}

Thus, the operator $L$ is well-defined.
$\square$

For the operator $L$, we now verify the following properties (i) - (iii):

(i) $L$ is a compact operator from $[\LOP, \, \UPP] \subset 
L^2(0,T;L^2(\OOO))$ into itself.  In other words, the set 
$L\, [\LOP, \, \UPP]$ is relatively compact in $L^2(0,T;L^2(\OOO))$.
\\
{\bf Verification of (i).}
Let $u \in [\LOP, \, \UPP]$.  Then the estimate \eqref{(6.12)} means 
$$
\Vert Lu\Vert_{L^2(0,T;H^2(\OOO))} 
+ \Vert Lu - a\Vert_{\HH(0,T;L^2(\OOO))} 
\le C_4.
$$
For $a\in H^1(\OOO)$, we deduce that
$$
\Vert Lu-a\Vert_{L^2(0,T;H^1(\OOO))}
+ \Vert Lu - a \Vert_{\HH(0,T;L^2(\OOO))}
\le C_4 + \sqrt{T}\Vert a\Vert_{H^1(\OOO)}.
$$
Since  
$$
L^2(0,T;H^1(\OOO)) \cap \HH(0,T;L^2(\OOO))
\subset L^2(0,T;L^2(\OOO))
$$ 
is a compact embedding (see, e.g., \cite{Te}), we obtain that 
$v \longrightarrow Lv - a$ is compact from $[\UPP,\, \LOP]
\subset L^2(0,T;L^2(\OOO))$ to $L^2(0,T;L^2(\OOO))$.
Thus the compactness of the operator $L$ is proved.
$\square$

(ii) $Lv \ge Lw$ in $\QQQQ$ if $v\ge w$ in $\QQQQ$ and
$v,w \in [\LOP, \,\UPP]$.
\\
{\bf Verification of (ii).}
Setting $y:= v-w$ and $z:= Lv - Lw$  
and 
applying the mean value theorem for  $f(v) - f(w)$, we obtain the 
representation
$$
\pppa (z-a) + Az + (M+1)z=(M+1)y 
+ f'(\mu(x,t))y \quad \mbox{in } \QQQQ,
$$
where $\mu(x,t)$ is a number between $v(x,t)$ and $u(x,t)$.

The inclusion  $f \in \mathcal{F}_M$ implies $\vert f'(\mu(x,t))\vert \le M$.
Therefore, in view of $y \ge 0$ in $\QQQQ$, we estimate
$$
(M+1)y + f'(\mu(x,t))y \ge (M+1-M)y(x,t) \ge 0 \quad \mbox{in } \QQQQ.
$$
Hence, the function $z$ satisfies
$$
\left\{ \begin{array}{rl}
& \pppa (z-a) + Az + (M+1)z \ge 0 \quad \mbox{in } \QQQQ,\\
& \NUNU z + \sigma z = 0 \quad \mbox{on } \ppp\OOO \times (0,T).
\end{array}\right.
$$
Theorem 1 from \cite{LYP1} implies that $z\in L^2(0,T;H^2(\OOO))$ and 
$z-a \in \HH(0,T;L^2(\OOO))$ and Theorem 2 from \cite{LYP1} yields the 
inequality $z\ge 0$ in $\QQQQ$.  Therefore $Lv \ge Lw$ in $\QQQQ$.
$\square$

(iii) $L[\LOP,\,\UPP] \subset [\LOP,\, \UPP]$.
\\
{\bf Verification of (iii).}
In order to prove the above inclusion,  we first show that $\LOP \le L\LOP$ .

Indeed, setting $v:= L\LOP$, we have
$$
\left\{ \begin{array}{rl}
& \pppa (v - a) 
+ Av + (M+1)v = (M+1)\LOP + f(\LOP) \quad \mbox{in } \QQQQ,\\
& v - a \in \HH(0,T;L^2(\OOO)), \\
& \NUNU v + \sigma u = 0 \quad \mbox{on } \ppp\OOO \times (0,T).
\end{array}\right.
$$
On the other hand, since $\LOP$ is a lower solution, it satisfies the following
relations:
$$
\left\{ \begin{array}{rl}
& \pppa (\LOP - \underline{a}) + A\LOP + (M+1)\LOP 
\le (M+1)\LOP + f(\LOP) \quad \mbox{in } \QQQQ,\\
& \LOP - \underline{a} \in \HH(0,T;L^2(\OOO)), \\
& \NUNU \LOP + \sigma \LOP = 0 \quad \mbox{on } \ppp\OOO \times (0,T).
\end{array}\right.
$$
Therefore, the function $z:= v - \LOP = L\LOP - \LOP$ is a solution to the 
problem
$$
\left\{ \begin{array}{rl}
& \pppa (z - (a -\underline{a})) + Az + (M+1)z \ge 0 
\quad \mbox{in } \QQQQ,\\
& z - (a-\underline{a}) \in \HH(0,T;L^2(\OOO)), \\
& \NUNU z + \sigma z = 0 \quad \mbox{on } \ppp\OOO \times (0,T),\\
& a-\underline{a} \ge 0 \quad \mbox{in } \OOO.
\end{array}\right.
$$
Then, according to Theorem 2 from \cite{LYP1}, the inequality 
$z\ge 0$ holds true in $\QQQQ$, that is, $L\LOP \ge \LOP$ in $\QQQQ$.
By similar arguments, we prove that $L\UPP \le \UPP$ in $\QQQQ$.

Now, let $\LOP \le u \le \UPP$ in 
$\QQQQ$.  
Since $L$ is increasing, we have the inequalities $L\LOP \le Lu \le L\UPP$.
Hence, $\LOP \le Lu \le \UPP$ in $\QQQQ$, which proves the property (iii).
$\square$

The properties (i)-(iii) allow to apply Lemmas \ref{l9.1} and \ref{l9.2},  
which completes the proof of 
Theorem \ref{t3.3}.
\end{proof}

Theorem \ref{t3.3} ensures that, for a given $T>0$, in order to guarantee the 
existence of the solution to the initial-boundary value problem \eqref{(2.3)}, 
it is sufficient 
to determine its upper and lower solutions.
This technique is called monotonicity method. It is related to the 
Perron method and to the concept of viscosity solutions.
For applications of the monotonicity method to the partial differential 
equations of parabolic type, 
that is, to the problem \eqref{(2.3)} with $\alpha=1$, we refer, e.g., to
\cite{Am}, \cite{Ke},  \cite{Pao1}, and \cite{Pao2}.  Theorem \ref{t3.3} 
asserts a corresponding result for the case 
$0<\alpha<1$. In the rest of this section, we discuss some important 
consequences that follow from this theorem. 

\begin{proposition}
\label{p3.1}
Let $u$ be a solution to the problem \eqref{(2.3)} satisfying  
\eqref{(3.4)} and the estimate $\vert u(x,t)\vert \le m$ for $(x,t)
\in \OOO\times (0,T)$ with a constant $m>0$. 

Then the following properties hold true:

(i) For any  upper solution $\UPP$ to the problem \eqref{(2.3)} satisfying
\eqref{(3.4)} and the estimate $\vert \UPP(x,t)\vert \le m$ for 
$(x,t) \in \OOO\times (0,T)$, 
we have
$$
u(x,t) \le \UPP(x,t) \quad \mbox{for all $(x,t)\in \OOO\times (0,T)$}.
$$ 

(ii) For any lower solution $\LOP$ to the problem \eqref{(2.3)} 
satisfying \eqref{(3.4)} and the estimate $\vert \LOP(x,t)\vert \le m$ for 
$(x,t) \in \OOO\times (0,T)$, we have
$$
\LOP(x,t) \le u(x,t) \quad \mbox{for all $(x,t)\in \OOO\times (0,T)$}.
$$ 
\end{proposition}

Henceforth we focus on a special case of the operator $-A$ from the 
formulation of the initial-boundary value problem \eqref{(2.3)}
in the form
\begin{equation}
\label{div}
\left\{\begin{array}{rl}
& -Av := \mbox{div}\, (p(x)\nabla v(x)) + c(x)v, \quad x\in \OOO,\\
& \DDD(A) = \{ v\in H^2(\OOO); \, \NUNU v + \sigma v = 0 \quad 
\mbox{on } \ppp\OOO \},
\end{array}\right.
\end{equation}
where $p\in C^2(\ooo{\OOO})$, $p(x)>0$ for $x\in \ooo{\OOO}$ and 
$c\in C^2(\ooo{\OOO})$, $c(x)\le 0$, $c(x) \not\equiv 0$ for $x \in \OOO$.


For this operator, Theorem \ref{t3.3} leads to the following result:

\begin{proposition}
\label{p3.2}
We assume that $f \in C^1(\R)$ is 
a monotone decreasing function and that 
there exists a solution $u_{\infty}\in H^2(\OOO)$ to the 
boundary-value problem
\begin{equation}
\label{(3.12)}
\left\{ \begin{array}{rl}
& Au_{\infty}(x) = f(u_{\infty}(x)), \quad x \in \OOO, \\
& \NUNU u_{\infty} + \sigma u_{\infty} = 0 \quad \mbox{on } \ppp\OOO,\\
& \vert u_{\infty}(x)\vert \le m, \quad x\in \OOO,
\end{array}\right.
\end{equation}
where the operator $-A$ is defined by \eqref{div}.

Then for any $T>0$ and any $a\in \DDD(A_0^{\gamma})$, 
the initial-boundary value problem \eqref{(2.3)} possesses
a unique solution $u_a$ in the class \eqref{(3.4)} such that the inequality 
\begin{equation}
\label{(3.14)}
\vert u_a(x,t) - u_{\infty}(x)\vert \le CE_{\alpha,1}(-\la_1 t^{\alpha})
\vert \va_1(x)\vert,
\quad x\in \OOO, \, t\in (0,T)                        
\end{equation}
holds true with  a certain constant $C>0$. 
In particular, we have the estimate
\begin{equation}
\label{(3.15)}
\vert u_a(x,t) - u_{\infty}(x)\vert \le Ct^{-\alpha}\vert \va_1(x)\vert
\quad \mbox{for all } x\in \OOO \mbox{ as } t\to \infty.   
\end{equation}
\end{proposition}

\begin{proof}

The proof is similar to the one of 
Proposition 1.4 (pp. 25-26) in \cite{Ke}, see also 
\cite{Pao2}.

Since the function $c\in C^2(\ooo{\OOO})$ satisfies $c(x)\le 0$ and 
$c(x)\not\equiv 0$ on $\ooo{\OOO}$, all eigenvalues of the operator 
$A$ are known to be positive and the sequence of the eigenvalues tends to 
$+\infty$: 
$$
0 < \la_1 \le \la_2 \le \cdots \longrightarrow \infty,
$$
where $\la_k$, $k\in \N$ are numbered according to their multiplicities.

Now we consider an eigenfuction $\va_1$ for the smallest  eigenvalue $\la_1$ 
satisfying the normalization condition $\Vert \va_1\Vert = 1$.  
It is well-known that 
the multiplicity of the eigenvalue $\la_1$ is one and
$\va_1(x) > 0$ for $x\in \ooo{\OOO}$ or
$\va_1(x) < 0$ for $x\in \ooo{\OOO}$ (see, e.g., 
Lemma 1.4 (p. 96) and
Theorem 1.2 (p. 97) in \cite{Pao2}).

For our proof,  we choose the positive eigenfunction $\va_1$ ($\va_1(x)>0$ 
for $x \in \ooo{\OOO}$).
Then we can find a sufficiently large constant $M_1>0$ such that 
\begin{equation}
\label{(9.2)}
u_{\infty}(x) - M_1\va_1(x) \le 
a(x) \le u_{\infty}(x) + M_1\va_1(x), \quad x\in \ooo{\OOO}.
\end{equation}
Now we  set 
$$
\left\{ \begin{array}{rl}
&\ooo{a}(x):= \UPP(x,0) = u_{\infty}(x) + M_1\va_1(x), \\
&\underline{a}(x):= u_{\infty}(x) - M_1\va_1(x), \quad x\in \OOO
\end{array}\right.
$$
and
$$
\left\{ \begin{array}{rl}
& \LOP(x,t):= u_{\infty}(x) - M_1\MLONE(-\la_1 t^{\alpha})\va_1(x),  \\
& \UPP(x,t):= u_{\infty}(x) + M_1\MLONE(-\la_1 t^{\alpha})\va_1(x),
\quad (x,t) \in \QQQQ,
\end{array}\right.
$$
where the Mittag-Leffler function $\MLONE$ is defined by the convergent series
$$
\MLONE(z) = \sum_{k=0}^{\infty} \frac{z^k}{\Gamma(\alpha k + 1)}, \quad 
\alpha >0,\ z\in \C.
$$
Let $T>0$ be arbitrary.  Recalling that 
$$
-Aw(x) = \mbox{div}\, (p(x)\nabla w) + c(x)w, \quad x\in \OOO
$$
with $c(x) \le 0$, $c(x) \not\equiv 0$ on $\ooo{\OOO}$ and 
$\pppa(\MLONE(-\la_1t^{\alpha})-1) = -\la_1\MLONE(-\la_1t^{\alpha})$,
we have the inclusion $\LOP - \lowa \in H_{\alpha}(0,T;L^2(\OOO))$ and 
the relations  
$$
\pppa (\LOP - \lowa) =  \la_1M_1\MLONE(-\la_1t^{\alpha})\va_1(x),
$$
$$
A\LOP = Au_{\infty} - \la_1M_1\MLONE(-\la_1t^{\alpha})\va_1(x),
\quad x\in \OOO,\, 0<t<T,
$$
so that 
$$
\pppa (\LOP - \lowa) + A\LOP = Au_{\infty} = f(u_{\infty}) \quad
\mbox{in }\QQQQ.
$$
It is known that $\MLONE(-\la_1t^{\alpha}) \ge 0$ for $t\ge 0$
(see, e.g., \cite{GKMR}) and 
$$
\LOP(x,t) \le u_{\infty}(x), \quad (x,t) \in \OOO\times (0,T).
$$
Since $f$ is decreasing, we have the inequality 
$f(u_{\infty}(x)) \le f(\LOP(x,t))$.  Therefore,  
$$
\pppa (\LOP - \lowa) + A\LOP \le f(\LOP) \quad
\mbox{in }\QQQQ.
$$
Since $\NUNU\LOP + \sigma\LOP = 0$ on 
$\ppp\OOO \times (0,T)$, the function 
$\LOP$ is a lower solution to the initial-boundary value problem \eqref{(2.3)} 
with the operator $-A$ defined by \eqref{div}. Similarly we can verify 
that $\UPP$ is an upper solution.  Hence, Theorem \ref{t3.3} yields that 
there exists a unique solution $u\in L^2(0,T;H^2(\OOO))$ such that 
$$
u-a \in \HH(0,T;L^2(\OOO)), \quad 
\LOP(x,t) \le u(x,t) \le \UPP(x,t), \quad (x,t) \in \QQQQ,
$$
that is,
$$
\vert u(x,t) - u_{\infty}(x)\vert \le M_1\MLONE(-\la_1t^{\alpha})
\va_1(x), \quad x\in \OOO,\, 0<t<T.
$$
Since $T>0$ is arbitrary and $\vert \MLONE(-\la_1t^{\alpha})\vert
\le \frac{C}{t^{\alpha}}$ as $t \to \infty$ (Theorem 1.6 (p. 35) in 
\cite{Po}), the proof of Proposition \ref{p3.2} is completed.
\end{proof}

In general, solution to the boundary-value problem \eqref{(3.12)} may  
be not unique. However,
under our assumptions,  the estimate \eqref{(3.14)} implies that $u_{\infty}$ 
is uniquely 
determined as the limit of $u_a$ as $t \to \infty$.
In the homogeneous  case ($f\equiv 0$), the asymptotic behavior 
of the solution $u_a$ is known  (see, e.g., \cite{KRY}, \cite{SY}, \cite{Ya}): 
$\Vert u_a(\cdot, t)\Vert = O(t^{-\alpha})$ as $t \to \infty$ that is 
the same as in the relation \eqref{(3.15)}.

We conclude this section with two examples of the results presented in 
Proposition \ref{p3.1} in the case
$A = -\Delta$ with the homogeneous Neumann boundary condition, that is, with 
$\sigma(x)=0$ on $\ppp\OOO$. 
In this case, the boundary condition can be represented as follows:
$$
\NUNU v = \ppp_{\nu}v := \nabla v\cdot \nu  = 0 \quad \mbox{on } \ppp\OOO.
$$

In what follows, we assume that there exists a solution $u_a$ to the 
initial-boundary value 
problem \eqref{(2.3)} that belongs the class defined as in \eqref{(3.4)}. 
In the examples, the estimates of
$u(x,t) - a(x)$ are derived for small $t>0$. 

\begin{example}
\label{e3}
Let 
$$
f(\eta) = -\frac{\eta}{1+\vert \eta\vert}, \quad \eta \in \R
$$
and the following conditions be satisfied:
\begin{equation}
\label{(3.16)}
a\in C^2(\ooo{\OOO}), \quad a(x)\ge 0 \quad \mbox{in $\OOO$}, \quad
\ppp_{\nu} a = 0 \quad\mbox{on }\ppp\OOO. 
\end{equation}

If the inequality 
\begin{equation}\label{(6.18a)}
\Delta a(x) \le \Gamma(\alpha+1)\rho, \quad x\in \ooo{\OOO}
\end{equation}
holds true for a constant $\rho>0$, then the estimates
\begin{equation}
\label{(3.17)}
0 \le u_a(x,t) - a(x) \le \rho t^{\alpha}, \quad x\in \OOO,\,
0\le t < T                         
\end{equation}
are satisfied. 

To prove the estimates \eqref{(3.17)}, 
we first set 
\begin{equation}
\label{(lu)}
\LOP(x,t) = \underline{a}(x) = 0, \quad
\UPP(x,t) = \rho t^{\alpha} + a(x), \quad 
x\in \OOO, \, 0<t<T.
\end{equation}
Then the equation 
$$
\pppa (\LOP - \underline{a}) - \Delta \LOP - f(\LOP) = 0
$$
evidently holds. Thus, the function $\LOP(x,t) \equiv 0$ satisfies the 
conditions
$$
\left\{ \begin{array}{rl}
& \pppa (\LOP - 0) - \Delta \LOP \le f(\LOP)\quad \mbox{in } \OOO
\times (0,T), \\
& 0=\underline{a}(x) \le u(x,0) = a(x), \quad x\in \OOO,\\
& \ppp_{\nu}\LOP = 0 \quad \mbox{on } \ppp\OOO \times (0,T),
\end{array}\right.
$$
which means that it is a lower solution to the initial-boundary value 
problem \eqref{(2.3)} with the data specified above. 

Moreover, the conditions \eqref{(3.16)} and \ref{(6.18a)} imply that the 
function $\UPP$ defined by \eqref{(lu)} is a solution to the problem
$$
\left\{ \begin{array}{rl}
& \pppa (\UPP(x,t) - a(x)) - \Delta \UPP - f(\UPP)
= \rho\Gamma(\alpha+1) - \Delta a 
+ \frac{a(x)+\rho t^{\alpha}}{1+\vert a(x)+\rho t^{\alpha}\vert} \\
\ge &\rho \Gamma(\alpha+1) - \Delta a(x) \ge 0, \quad x\in \OOO, \, 0<t<T, \\
& \UPP(x,0) = a(x), \quad x\in \OOO,\\
& \ppp_{\nu} \UPP = 0\quad \mbox{on } \ppp\OOO \times (0,T). 
\end{array}\right.
$$
Therefore, the functions $\UPP$ and $\LOP$ defined by \eqref{(lu)} are 
an upper and a lower solutions,
respectively, and Proposition \ref{p3.1} immediately leads to the estimate 
\eqref{(3.17)}.
\end{example}

We note that for $\alpha=1$ the equation considered in Example \ref{e3} is a
semilinear equation of parabolic type
$$
\ppp_tu - \Delta u = -\frac{u}{1+\vert u\vert},
$$
which is known to be a governing equation for an enzyme process. For this type 
of semilinear equations of parabolic type, the monotonicity method applied above is well-known  
(see, e.g., \cite{Ke} or \cite{Pao1}).

\begin{example}
\label{e4}
Let $f \in C^1(0,\infty)$  be an increasing function 
in $\R_{+}$ and
\begin{equation}
\label{(3.18)}
a \in C^2(\ooo{\OOO}), \qquad \ppp_{\nu}a = 0 \quad \mbox{on } \ppp\OOO.  
\end{equation}
Then, for an arbitrary but fixed $\ep \in (0, \alpha)$, there exists 
a positive number $T_1=T_1(\ep) > 0$ such that 
\begin{equation}
\label{(3.19)}
u_a(x,t) - a(x) \le t^{\alpha-\ep}, \quad 0<t<T_1.     
\end{equation}
\end{example}

In order to prove the inequality \eqref{(3.19)}, 
we set $\UPP(x,t) = t^{\alpha-\ep} + a(x)$ for $(x,t) \in 
\OOO\times (0,T_1)$ with $T_1$ which will be selected later.  Then 
$$
\pppa (\UPP-a) = \frac{\Gamma(\alpha-\ep+1)}{\Gamma(-\ep+1)}t^{-\ep}
$$
and $\Delta \UPP = \Delta a$, and thus
$$
\pppa (\UPP - a) - \Delta \UPP 
= \frac{\Gamma(\alpha-\ep+1)}{\Gamma(-\ep+1)}t^{-\ep} - \Delta a,
\quad x\in \OOO, \, 0<t<T_1.
$$
Moreover, setting $M_1:= \max_{x\in \ooo{\OOO}} a(x)$, we get
\begin{equation}\label{(6.23)}
f(\UPP(x,t)) = f(t^{\alpha-\ep}+a(x))
\le f(T_1^{\alpha-\ep}+M_1), \quad x\in \OOO,\, 0<t<T_1.
\end{equation}
We note that $T_1^{-\ep}$ can be arbitrarily large if we choose 
$T_1>0$ sufficiently small.  Therefore, because $\alpha-\ep>0$, 
we can choose $T_1 > 0$ sufficiently small 
such that 
\begin{equation}
\label{(3.20)}
\frac{\Gamma(\alpha-\ep+1)}{\Gamma(-\ep+1)}T_1^{-\ep}
\ge f(T_1^{\alpha-\ep} + M_1) + \Delta a(x), \quad x\in \ooo{\OOO}.
\end{equation}
For this $T_1$, the inequalities \eqref{(6.23)} and \eqref{(3.20)} readily 
lead to the estimate
$$
\pppa (\UPP-a) - \Delta \UPP \ge f(\UPP) \quad \mbox{in } \OOO\times (0,T_1).
$$
Since $\UPP(x,0) = a(x) = u_a(x,0)$ and $\ppp_{\nu}\UPP = 0$ on 
$\ppp\OOO \times (0,T_1)$, the function  $\UPP$ is an upper solution for
$t\in (0,\, T_1)$ that immediately yields the inequality \eqref{(3.19)}.

Next, we look for a suitable  lower solution. 
In addition to the conditions \eqref{(3.18)}, we assume that there exists
a constant $\delta_1>0$ such that 
\begin{equation}\label{(6.23a)}
a(x) \ge \delta_1 > 0 \quad \mbox{for all $x\in \ooo{\OOO}$}.
\end{equation}
Then we set  
\begin{equation}
\label{(3.21)}
M_2 := \max_{x\in \ooo{\OOO}} (-\Delta a(x)).     
\end{equation}
The lower solution is introduced in the form
$$
\LOP(x,t) = a(x) - \rho t^{\alpha}, \quad x \in \OOO,\,
0<t < \left( \frac{\delta_1}{2\rho}\right)^{\frac{1}{\alpha}}=: T_2,
$$
where $\delta_1>0$ is the constant from the condition \eqref{(6.23a)} and 
the constant $\rho>0$ will be selected later.
Then, \eqref{(6.23a)} implies
$$
\LOP(x,t) \ge \frac{\delta_1}{2}, \quad x\in \OOO,\, 0<t<T_2.
$$
Moreover, the inequalities 
$$
\pppa (\LOP - a) - \Delta \LOP = -\Gamma(\alpha+1)\rho - \Delta a
\le -\Gamma(\alpha+1)\rho + M_2
$$
and
$$
f(\LOP(x,t)) = f(a(x) - \rho t^{\alpha})
\ge f(\delta_1 - \rho T_2^{\alpha}) 
= f\left( \delta_1 - \frac{\delta_1}{2}\right) 
= f\left(\frac{\delta_1}{2}\right)
$$
hold true.
Now we choose a constant $\rho$ sufficiently large such that 
$T_2:= \left( \frac{\delta_1}{2\rho}\right)^{\frac{1}{\alpha}}$
is sufficiently small and the inequality
$$
\rho \ge \frac{1}{\Gamma(\alpha+1)}
\left(M_2 - f\left(\frac{\delta_1}{2}\right)\right)
$$
holds true. 
Then we have the estimate
$$
\pppa (\LOP - a) - \Delta \LOP \le f(\LOP) \quad \mbox{in } \OOO\times 
(0,T_2).
$$
Since $\ppp_{\nu}\LOP = 0$ on $\ppp\OOO \times (0,T_2)$ and 
$\LOP(x,0) = a(x)$ for $x\in \OOO$, the function
$\LOP$ is a lower solution, which leads to the inequality
\begin{equation}
\label{(3.22)}
a(x) - \frac{1}{\Gamma(\alpha+1)}
\left(M_2 - f\left(\frac{\delta_1}{2}\right)\right)t^{\alpha} \le u_a(x,t) 
\end{equation}
for $0 < t < T_2$, where $T_2>0$ is sufficiently small.

Summarizing the results presented above,  under the conditions  \eqref{(3.18)} 
and \eqref{(6.23a)}, for any $\ep \in (0,\, \alpha)$, there exists 
a constant $T_{\ep,a} > 0$ such that  
\begin{equation}
\label{(3.23)}
-\frac{1}{\Gamma(\alpha+1)}\left( M_2 - f\left( \frac{\delta_1}{2}\right)
\right)t^{\alpha} \le u_a(x,t) - a(x)
\end{equation}
$$
\le t^{\alpha-\ep}, \quad x\in \OOO,\, 0<t<T_{\ep,a}.                     
$$

Finally, let us emphasize that the constructions of the upper and lower solutions 
presented above  are not unique.  For example,
the right-hand side of the inequality \eqref{(3.23)} has the form 
$t^{\alpha-\ep}$. However, under the conditions \eqref{(3.18)}, we can obtain 
a different estimate as follows: For a sufficiently small $T_3>0$, 
we can choose 
a constant $M_3=M_3(T_3)>0$ such that 
\begin{equation}\label{(6.27a)}
\left\{ \begin{array}{rl}
&\Vert \Delta a\Vert_{C(\ooo{\OOO})} \le \frac{1}{2}M_3\Gamma(\alpha+1),\\
& f(M_3T_3^{\alpha}+\Vert a\Vert_{C(\ooo{\OOO})})
\le \frac{1}{2}M_3\Gamma(\alpha+1).
\end{array}\right.
\end{equation}
Indeed, for a small $T>0$, the relation 
$\lim_{T\downarrow 0} f(M_3T^{\alpha}+\Vert a\Vert_{C(\ooo{\OOO})})
= f(\Vert a\Vert_{C(\ooo{\OOO})})$ ensures existence of a sufficiently large 
constant $M_3(T)>0$ that fulfills the inequalities \eqref{(6.27a)}.

Then the inequality 
\begin{equation}
\label{(3.24)}
u_a(x,t) - a(x) \le M_3(T_3)t^{\alpha},\quad x\in \OOO, \, 0<t<T_3
\end{equation}
holds true. 

To deduce it, 
we set $\UPP(x,t) := a(x) + M_3t^{\alpha}$.  Then the first condition 
in \eqref{(6.27a)} yields 
$$
\pppa (\UPP-a) - \Delta \UPP = M_3\Gamma(\alpha+1) - \Delta a
\ge \frac{1}{2}M_3\Gamma(\alpha+1).
$$
Moreover, since $f$ is monotone increasing, we have the estimate
$$
f(\UPP) = f(a(x)+M_3t^{\alpha}) 
\le f(\Vert a\Vert_{C(\ooo{\OOO})}+ M_3T_3^{\alpha})
$$
for $(x,t) \in\OOO \times (0,T_3)$.
Therefore, by means of \eqref{(6.27a)}, we obtain
$$
\pppa (\UPP-a) - \Delta \UPP \ge \frac{1}{2}M_3\Gamma(\alpha+1)
\ge f(M_3T_3^{\alpha}+\Vert a\Vert_{C(\ooo{\OOO})})
\ge f(\UPP) \quad \mbox{in $\OOO \times (0,T_3)$},
$$
which means that $\UPP$ is another upper solution and 
thus the inequality \eqref{(3.24)} is verified.

\section{Non-negativity results for systems of the linear time-fractional diffusion equations}

\label{sec:5}

\setcounter{section}{5} 
\setcounter{equation}{0} 

In this section, we treat the initial-boundary value problems for the systems of the linear time-fractional diffusion equations in the form
$$
\left\{ \begin{array}{rl}
& \left(
\begin{array}{c}
\ppp_t^{\alpha_1} (u_1(x,t) - a_1(x))  \\
 \cdots                      \\
 \ppp_t^{\alpha_N} (u_N(x,t) - a_N(x))  \\
\end{array}\right)
= \left(\begin{array}{c}
\Delta u_1(x,t)  \\
 \cdots \cdots       \\
\Delta u_N(x,t) \\
\end{array}\right)              \cr\\            
+ &\left( \begin{array}{ccc}
 p_{11}(x,t) &\cdots &  p_{1N}(x,t) \\
      \cdots &\cdots &  \cdots      \\
 p_{N1}(x,t) & \cdots & p_{NN}(x,t)  \\
\end{array}\right)
\left(
\begin{array}{c}
u_1(x,t) \\
\cdots                      \\
u_N(x,t) \\
\end{array}\right)             
+ \left(
\begin{array}{c}
F_1(x,t)   \\
 \cdots                      \\
F_N(x,t) \\
\end{array}\right),             \cr \\
& u_k - a_k \in H_{\alpha_k}(0,T;L^2(\OOO)), \cr\\
& \ppp_{\nu}u_k = 0 \qquad \mbox{on $\ppp\OOO \times (0,T)$ for 
$1\le k \le N$,}
\end{array}\right.                                            \eqno{(5.1)}
$$
\noindent
where $N\ge 2$ and the orders $\alpha_n,\ n=1,\dots N$ of the fractional derivatives from the left-hand side of the equations (5.1) satisfy the conditions
$$
0<\alpha_1 < \alpha_2 < \cdots < \alpha_N < 1.
$$

For the sake of simplicity, in the equations (5.1), we restrict ourselves to the case of the spatial differential operator  $\Delta := 
\sum_{j=1}^d \ppp_j^2$  coupled with the zeroth order terms.
The case of the time-fractional linear diffusion equations with the homogeneous Dirichlet boundary condition was discussed in \cite{LY1}.
For $a_k\in H^1(\OOO)$ and $F_k \in L^2(\QQQQ)$, $1\le k \le N$,
existence of a unique solution $(u_1, ..., u_N) \in 
(L^2(0,T;H^2(\OOO))^N$ to the equations (5.1) such that
$(u_1-a_1, \, \cdots u_N-a_N) \in (H_{\alpha_1}(0,T;L^2(\OOO)) \times
\cdots \times H_{\alpha_N}(0,T;L^2(\OOO)))$ can be proven in analogy to the derivations presented in \cite{LY1} and \cite{LYP1} for the case of the time-fractional linear diffusion equations with the homogeneous Dirichlet boundary condition. 
  
The main result of this section is formulated below.

\begin{theorem}

For $1\le j,k \le N$, let $p_{jk} \in L^{\infty}(\QQQQ)$, and
$$
p_{jk} \ge 0 \quad \mbox{in $\OOO\times (0,T)$} \quad \mbox{if $j\ne k$},
                                         \eqno{(5.2)}
$$
$$
\left\{ \begin{array}{rl}
& F_k \in L^2(\OOO\times (0,T)), \quad F_k \ge 0 \quad
\mbox{in $\OOO\times (0,T)$}, \cr\\
& a_k \in H^1(\OOO), \quad a_k \ge 0 \quad \mbox{in $\OOO$}.
\end{array}\right.
                             \eqno{(5.3)}
$$
Then $u_k(x,t) \ge 0$ for $(x,t) \in \OOO\times (0,T)$ for 
$1\le k \le N$.
\end{theorem}

It is worth mentioning that the conditions (5.2) allow increasing of $\ppp_t^{\alpha_{\ell}}(u_\ell - a_\ell)$ if all components $u_k$ with $k \ne \ell$ increase. Thus
the system (5.1)  can be characterized as a cooperative one and then one can expect 
the non-negativity of all components $u_{\ell}$ if the
initial values are  non-negative.  In the case of the parabolic PDEs ($\alpha_1 = \alpha_2 = \cdots = \alpha_N = 1$), (5.1) is indeed known to be a cooperative system.
\\
{\bf Sketch of a proof of Theorem 5.1.}
\\
Henceforth we introduce the notations
$$
A_0v(x) := - \Delta v + M_1v, \quad 
\DDD(A_0) = \{v \in H^2(\OOO);\, \ppp_{\nu}v\vert_{\ppp\OOO} = 0 \}.
$$
In the above definition, we choose a constant $M_1 > 0$ such that 
$p_{jj}(x,t) + M_1 > 0,\  j=1,2,\dots,N$ for $(x,t)\in \QQQQ$. 
Due to the condition $M_1>0$, all  eigenvalues of the operator $A_0$ are positive and can be ordered according to their multiplicities:
$$
0 < \la_1 \le \la_2 \le \cdots \longrightarrow \infty.
$$
In what follows, by $\va_n,\ n=1,2,\dots$ we denote the linearly 
independent eigenfunctions of $A_0$ corresponding to the eigenvalues $\la_n$ such that 
$\{ \va_n\}_{n\in \N}$ is an orthonormal basis in $L^2(\OOO)$.

Then, in analogy to (3.8) and (3.9), we set 
$$
\left\{ \begin{array}{rl}
& S_{\ell}(t)a := \sumn E_{\alpha_{\ell},1}(-\la_nt^{\alpha_\ell})(a,\va_n)
\va_n,                             \cr\\
& K_{\ell}(t)a := \sumn t^{\alpha_{\ell}-1}
E_{\alpha_{\ell},\alpha_{\ell}}(-\la_nt^{\alpha_\ell})(a,\va_n)\va_n, \ 
t>0, \, a\in L^2(\OOO).
\end{array}\right.
                                \eqno{(5.4)}
$$
Noting that 
$t^{\alpha_\ell-1} \le Ct^{\alpha_1-1}$ for $\ell=1,..., N$ and acting as by derivation of the estimates (3.10), 
we can prove the inequalities
$$
\left\{
\begin{array}{rl}
& \Vert S_{\ell}(t)a\Vert \le C\Vert a\Vert,\ \ell=1,..., N, \cr\\
& \Vert K_{\ell}(t)a \Vert \le Ct^{\alpha_1-1}\Vert a\Vert,\ \ell=1,..., N \quad
\mbox{for all $a\in L^2(\OOO)$}.
\end{array}\right.
                                  \eqno{(5.5)}
$$
Based on these estimates, we  verify the inclusions 
$$
S_{\ell}(t)a \in L^\infty(0,T;L^2(\OOO)) \quad \mbox{for all
$a \in L^2(\OOO)$}
$$
and
$$
\int^t_0 K_{\ell}(t-s)F(s) ds \in L^2(\QQQQ) 
\quad \mbox{for all $F\in L^2(\QQQQ)$}.
$$
In particular, the second inclusion follows from the Young inequality 
for the Laplace convolution.

By Theorem 2 in \cite{LYP1}, for each $\ell=1, ..., N$, we arrive at the inequalities
$$
S_{\ell}(t)a \ge 0 \quad \mbox{in $\QQQQ$ if $a \in H^1(\OOO)$ and 
$a \ge 0$ in $\OOO$}              \eqno{(5.6)}
$$
and
$$
\int^t_0 K_{\ell}(t-s)F(s) ds \ge 0 \quad \mbox{in $\QQQQ$}
$$
$$
\mbox{if 
$F\in L^2(\QQQQ)$ and $F\ge 0$ in $\QQQQ$.}
                                              \eqno{(5.7)}
$$
Indeed, the functions $v_1(t):= S_{\ell}(t)a$ and $v_2(t):= \int^t_0 K_{\ell}
(t-s) F(s) ds$ are solutions to the boundary-value problems
$$
\left\{ \begin{array}{rl}
& \ppp_t^{\alpha_{\ell}} (v_1 -a) = -A_0v_1 \quad \mbox{in $\QQQQ$}, \cr\\
& \ppp_{\nu}v_1 = 0 \quad \mbox{on $\ppp\OOO \times (0,T)$}
\end{array}\right.
$$
and
$$
\left\{ \begin{array}{rl}
& \ppp_t^{\alpha_{\ell}} v_2 = -A_0v_2 + F \quad \mbox{in $\QQQQ$}, \cr\\
& \ppp_{\nu}v_1 = 0 \quad \mbox{on $\ppp\OOO \times (0,T)$}
\end{array}\right.
$$
respectively. The Theorem 2 in \cite{LYP1} yields then the estimates (5.6) and (5.7).

Now we construct a sequence $U^n=(u_1^n, ..., u^n_N)$, $n=0,1,2,\dots$ as follows:
$$
\left\{ \begin{array}{rl}
& U^0(x,t)= (u_1^0(x,t), ..., u_N^0(x,t)) := (a_1(x), ..., a_n(x)) 
\quad \mbox{in $\OOO\times (0,T)$},                       \cr\\
& \ppp_t^{\alpha_1}(u_1^n - a_1) - \Delta u_1^n 
+ M_1u_1^n = (p_{11}+M_1)u_1^{n-1} + \cdots + p_{1N}u_{N}^{n-1} + F_1,
                                            \cr\\
& \cdots \cdots                                    \cr\\
& \ppp_t^{\alpha_N}(u_N^n - a_N) - \Delta u_N^n 
+ M_1u_N^n = p_{N1}u_1^{n-1} + \cdots + p_{N,N-1}u_{N-1}^{n-1} 
                                        \cr\\
+ &(p_{NN} + M_1)u_N^{n-1} + F_N, \qquad \mbox{in $\OOO\times (0,T)$}, \cr\\
& \ppp_{\nu}u_\ell^n = 0 \quad \mbox{on $\ppp\OOO \times (0,T)$, 
$1\le \ell \le N$}. 
\end{array}\right.
$$
Henceforth we employ the notation $u(t):= u(\cdot,t)$ and interpret it as
a function  mapping $t$ to an element from $L^2(\OOO)$.

Then, in view of the relations  (2.14) and (4.3)  from \cite{LYP1}, 
we can represent $u_\ell^n$, $n\in \N$ in terms of the operators $S_{\ell}(t)$ and $K_{\ell}(t)$: 
$$
\left\{ \begin{array}{rl}
& u_\ell^n(t) = S_\ell(t)a_\ell 
+ \int^t_0 K_\ell(t-s)\left( \sum_{j=1}^N (p_{\ell j}(s)
+ \delta_{\ell j}M_1)u_j^{n-1}(s) \right) ds   \cr\\\
+ & \int^t_0 K_\ell(t-s)F_{\ell}(s) ds \quad \mbox{in 
$0<t<T$ for $1\le \ell \le N$},\cr\\
& (u_1^0(t), ..., u^0_N(t)) = (a_1, ..., a_N).
\end{array}\right.
                                                  \eqno{(5.8)}
$$
Here and henceforth we set $\delta_{\ell j} = 0$ if $\ell \ne j$ and
$\delta_{jj} = 1$.

The representations (5.8) and the estimates (5.5) ensure that the sequence $U^n = (u_1^n, ..., u_N^n)
\in (L^2(\QQQQ))^N$ is well-defined for each 
$n\in \N$.

Moreover, starting from the estimates $a_{\ell} \ge 0,\ \ell=1, ..., N$ in $\QQQQ$ and employing the inequalities 
$p_{\ell j} + \delta_{\ell j}M_1 \ge 0,\ \ell, j = 1, ..., N$ in $\QQQQ$,  by the principle of mathematical induction, we  deduce the non-negativity of all components of $U^n,\ n=1, ..., N$:
$$
u_1^n, \, ..., u_N^n \ge 0 \quad \mbox{in $\QQQQ$}.   \eqno{(5.9)}
$$

Next we prove convergence of the sequence $U^n,\ n=0,1,2,\dots$ in $(L^2(\QQQQ))^N$, i.e., 
$$
U^n = (u_1^n, ..., u_N^n) \longrightarrow (u_1, ..., u_N) \quad
\mbox{in $(L^2(\QQQQ))^N$ as $n\to \infty$}.
                                     \eqno{(5.10)}
$$

Once the relation (5.10) is proved, the proof of Theorem 5.1 is completed. Indeed,
because of (5.10), the components of $(u_1, ..., u_N)$ are solutions to the equations
\begin{align*}
& u_\ell(t) = S_\ell(t)a_\ell 
+ \int^t_0 K_\ell(t-s)\left( \sum_{j=1}^N (p_{\ell j}(s)
+ \delta_{\ell j}M_1)u_j(s) \right) ds   \cr\\\
+ & \int^t_0 K_\ell(t-s)F_{\ell}(s) ds \quad \mbox{in 
$0<t<T$ for $1\le \ell \le N$}.
\end{align*}

Moreover, we can choose a subsequence $n' $ in  $\N$ such that 
$u_{\ell}^{n'} \longrightarrow u_{\ell}$ almost everywhere in $\QQQQ$ for 
each $\ell=1. ..., N$.  Hence the inequalities (5.9) imply
$u_{\ell} \ge 0,\ \ell=1. ..., N$ in $\QQQQ$ that completes the proof of Theorem 5.1.

Now let us prove the relation  (5.10).

Employing the equation (5.8), for $0<t<T$ and $\ell=1. ..., N$, we obtain the representation
$$
 u_\ell^{n+1}(t) - u_{\ell}^n(t) 
= \int^t_0 K_\ell(t-s)\left( \sum_{j=1}^N (p_{\ell j}(s) 
+ \delta_{\ell j}M_1)(u_j^n(s) - u_j^{n-1}(s))\right) ds.
$$
The inequalities (5.5) and the inclusions $p_{\ell j} \in L^{\infty}(\QQQQ)$  
lead then to the estimates
$$
 \Vert u_\ell^{n+1}(t) - u_{\ell}^n(t) \Vert
\le C\int^t_0 (t-s)^{\alpha_1-1} \left( 
\sum_{j=1}^N \Vert u_j^{n}(s) - u_j^{n-1}(s))\Vert \right) ds,\ \ell=1. ..., N.
$$

Setting $U_n(t) := \sum_{j=1}^N \Vert u_j^{n+1}(t) 
- u_j^n(t)\Vert$ for $n\in \N$, we rewrite the last inequality in the form
$$
U_n(t) \le C(J^{\alpha_1}U_{n-1})(t), \quad 0<t<T,\, n\in \N,
                                           \eqno{(5.11)}
$$
where $J^{\alpha}$ is the Riemann-Liouville fractional integral defined by  $J^{\alpha}v(t) := \frac{1}{\Gamma(\alpha)}
\int^t_0 (t-s)^{\alpha-1}v(s) ds$. 

Now we proceed with estimations of $U_n$ in the same manner as in the proof of Theorem 1 in \cite{LYP1}.

Since $J^{\alpha_1}v_1(t) \le J^{\alpha_1}v_2(t)$ for $0\le t \le T$ if 
$v_1(t) \le v_2(t)$ for $0\le t \le T$, we get
$$
U_1(t) \le C(J^{\alpha_1}U_0)(t)
= \frac{CU_0}{\Gamma(\alpha_1)}\int^t_0 (t-s)^{\alpha_1-1}ds 
= \frac{CU_0}{\Gamma(\alpha_1+1)}t^{\alpha_1}
$$
and then
$$
U_2(t) \le C(J^{\alpha_1}U_1)(t) 
\le \frac{C}{\Gamma(\alpha_1)}\int^t_0 (t-s)^{\alpha_1-1}\cdot
\frac{CU_0}{\Gamma(\alpha_1+1)}s^{\alpha_1} ds 
= \frac{C^2U_0}{\Gamma(2\alpha_1+1)}t^{2\alpha_1}.
$$
Using the same procedure $n$ times, we arrive at the estimates
$$
U_n(t) \le \frac{C^nU_0}{\Gamma(n\alpha_1+1)}t^{n\alpha_1}
\le \frac{C^nU_0}{\Gamma(n\alpha_1+1)}T^{n\alpha_1}, \ 0\le t\le T
\eqno{(5.12)}
$$
valid for any $n=1,2,3,\dots$.

Due to the 
 known asymptotic formula 
$$
\frac{\Gamma(x+a)}{\Gamma(x+b)} \sim x^{a-b},\ x \to +\infty
$$
for the Gamma function $\Gamma$, we then get the relation
$$
\lim_{n\to+\infty} 
\frac{ \frac{C^{n+1}T^{(n+1)\alpha_1}}{\Gamma((n+1)\alpha_1+1)} }
{   \frac{C^{n}T^{n\alpha_1}}{\Gamma(n\alpha_1+1)} } = \lim_{n\to+\infty}  \frac{C\, T^{\alpha_1}}{ (n\alpha_1)^{\alpha_1}}
= 0.
$$
Combining it with the estimates (5.12) we conclude that the series
$\sum_{n=1}^{\infty} U_n(t)$ is convergent in $L^{\infty}(0,T)$. Thus the relation  (5.10) is proved that completes  the proof of Theorem 5.1.
$\square$
\section{Non-negativity results for systems of the semilinear time-fractional diffusion equations}
\label{sec:6}

\setcounter{section}{6} 
\setcounter{equation}{0} 

In this section, we demonstrate how to apply Theorem 5.1 presented in the previous section  for analysis of the initial-boundary value problems for  the systems  of the  semilinear time-fractional diffusion equations. For the sake of simplicity, in what follows, we restrict ourselves to the following system of two equations even if  the same technique is applicable for the systems with  arbitrary many equations:  
$$
\left\{ \begin{array}{rl}
& \pppa (u-a) = \Delta u + f(u,v) \quad \mbox{in $\OOO\times (0,T)$}, \\
& \pppa (v-b) = \Delta v + g(u,v) \quad \mbox{in $\OOO\times (0,T)$},\\
& \ppp_{\nu}u = \ppp_{\nu}v = 0 \quad \mbox{on $\ppp\OOO\times (0,T)$}.
\end{array}\right.
                             \eqno{(6.1)}
$$
In (6.1), $\OOO \subset \R^3$ is a bounded domain with a smooth boundary 
$\ppp\OOO$ and the order $\alpha$ of the fractional derivatives satisfies the conditions $0<\alpha <1$. 

For formulation of further conditions, we use the notations introduced in Section 5: $A_0v := -\Delta v + M_1v$ for $v\in \DDD(A_0)
:= \{ w\in H^2(\OOO);\, \ppp_{\nu}w\vert_{\ppp\OOO} = 0\}$, where 
$M_1>0$ is a sufficiently large constant.

In what follows, we fix a constant $\gamma$ such that   
$
\frac{3}{4} < \gamma < 1,
$
and assume that the semilinear terms $f,g \in C^1(\R^2)$ satisfy 
the conditions that are  similar to the condition  (3.2) from Section 3: for some constant $m>0$, there exists
a constant $C = C(f,g,m) >0$ such that  
$$
\left\{ \begin{array}{rl}
& \mbox{(i)} \quad \Vert f(u_1,v_1) - f(u_2,v_2)\Vert
  + \Vert g(u_1,v_1) - g(u_2,v_2)\Vert  \cr\\
\le & C(\Vert u_1 - u_2\Vert_{\DDD(A_0^{\gamma})} 
     + (\Vert v_1 - v_2\Vert_{\DDD(A_0^{\gamma})}) \cr\\
\quad & \mbox{if $\Vert u_j\Vert_{\DDD(A_0^{\gamma})}, \, \,  
 \Vert v_j\Vert_{\DDD(A_0^{\gamma})} \le m$ for $j=1,2$}.\cr\\
& \mbox{(ii) There exists a constant $\ep>0$ such that}    \cr\\
& \Vert f(u,v)\Vert_{H^{2\ep}(\OOO)}, \, \,
 \Vert g(u,v)\Vert_{H^{2\ep}(\OOO)} \le C \cr\\
& \mbox{if $\Vert u\Vert_{\DDD(A_0^{\gamma})}$ and
           $\Vert v\Vert_{\DDD(A_0^{\gamma})} \le m$}.
\end{array}\right.
                                      \eqno{(6.2)}
$$

Then, for any $a, b$ satisfying the estimates $\Vert a\Vert_{\DDD(A_0^{\gamma})}\le m$,
$\Vert b\Vert_{\DDD(A_0^{\gamma})} \le m$ with a constant $m>0$,  there exists a number  $T=T(m)>0$ such that the problem 
(6.1) possesses a unique solution $(u,v) \in \{ L^2(0,T;H^2(\OOO)) \cap
C([0,T]; \DDD(A_0^{\gamma})) \}^2$ satisfying the initial conditions
$u-a, v-b \in H_{\alpha}(0,T;L^2(\OOO))$. The proof of this statement follows the lines of the proof of Theorem 3.4 in Section 3 and we omit it here. 

Now we formulate and prove the main result of this section.

\begin{theorem}
\label{t5}
Let $(u,v) \in \{ L^2(0,T;H^2(\OOO)) \cap C([0,T];\DDD(A_0^{\gamma}))
\}^2$ be a solution to the problem (6.1) 
such that 
$u-a, v-b \in H_{\alpha}(0,T;L^2(\OOO))$.

Then 
$u\ge 0$ and $v\ge 0$ in $\OOO\times (0,T)$ provided $a \ge 0$ and $b\ge 0$ in $\OOO$ and the conditions
$$
f(0,\eta) \ge 0 \quad \mbox{or}\quad
\left\{ \begin{array}{rl}
& (\ppp_2f)(\xi,\eta) \ge 0, \cr\\
& f(0,0) = 0
\end{array}\right.
\quad \mbox{for $\xi, \eta \in \R$}
                                   \eqno{(6.3)}
$$
and
$$
g(\xi,0) \ge 0 \quad \mbox{or}\quad
\left\{ \begin{array}{rl}
& (\ppp_1g)(\xi,\eta) \ge 0, \cr\\
& g(0,0) = 0
\end{array}\right.
\quad \mbox{for $\xi, \eta \in \R$}
                                   \eqno{(6.4)}
$$
are satisfied. 

\end{theorem}

\begin{proof}
For $\gamma > \frac{3}{4}$, the Sobolev embedding ensures the inclusion $\DDD(A_0^{\gamma}) \subset C(\ooo\OOO)$ and thus
we can choose a constant $M>0$ such that 
$$
\vert u(x,t)\vert, \, \vert v(x,t)\vert \le M \quad
\mbox{for $(x,t) \in \OOO\times (0,T)$}.
$$
Since $f,g \in C^1(\R^2)$, there exists a constant $C_1 = C_1(M) > 0$
such that 
$$
\vert \ppp_kf(\xi, \eta)\vert, \, \vert \ppp_kg(\xi,\eta)\vert 
\le C_1 \quad \mbox{for $k=1,2$ if $\xi,\eta \in [-M,M]^2$}.
                                 \eqno{(6.5)}
$$

Now we separately consider four possibilities for realization of the condition (6.3) and (6.4).
\\
{\bf Case 1:}
$$
f(0,\eta) \ge 0, \quad g(\xi,0) \ge 0 \quad 
\mbox{for $\xi, \eta \in \R$.}
$$
Then the mean value theorem yields
$$
f(u(x,t),v(x,t)) = f(u,v) - f(0,v) + f(0,v)
= (\ppp_1f)(\xi_1,v)u + f(0,v)        \eqno{(6.6)}
$$
and
$$
g(u(x,t),v(x,t)) = g(u,v) - g(u,0) + g(u,0)
= (\ppp_2g)(u, \eta_1)v + g(u,0),        \eqno{(6.7)}
$$
where $\xi_1 = \xi_1(u(x,t))$ is a number between $u(x,t)$ and $0$,
and $\eta_1 = \eta_1(v(x,t))$ is a number between $v(x,t)$ and $0$.

As a result, the problem (6.1) can be rewritten in the form
$$
\left\{ \begin{array}{rl}
& \pppa (u-a) = \Delta u + p_1(x,t)u + f(0,v) \ \ \mbox{in $\OOO\times (0,T)$}, \cr\\
& \pppa (v-b) = \Delta v + q_1(x,t)v + g(u,0)
\ \ \mbox{in $\OOO\times (0,T)$}, \cr\\
& \ppp_{\nu}u = \ppp_{\nu}v = 0 \quad \mbox{on $\ppp\OOO\times (0,T)$}, 
\end{array}\right.
$$
where we set 
$$
p_1(x,t) := (\ppp_1f)(\xi_1(u(x,t)),v(x,t)), \quad 
q_1(x,t) := (\ppp_2g)(u(x,t),\eta_1(v(x,t)).
$$
Due to the relation (6.5), the inclusions $p_1, q_1 \in L^{\infty}(\OOO\times (0,T))$ hold valid.
Since $f(0,v(x,t)) \ge 0$ and $g(u(x,t), 0) \ge 0$ for 
$(x,t)\in \OOO\times (0,T)$ by the assumptions made for Case 1, we can apply Theorem 5.1 and  get the non-negativity of the solutions to the problem (6.1): $u, v \ge 0$ in $\OOO\times (0,T)$.
\\
{\bf Case 2:}
$$
f(0,\eta) \ge 0, \quad \ppp_1g(\xi,\eta) \ge 0, \quad g(0,0) = 0 \quad
\mbox{for $\xi, \eta \in \R$.}
$$

In this case, the relation (6.6)  for the function $f$ holds valid. For the function $g$, we get the representation
$$
g(u(x,t), v(x,t)) = g(u,v) - g(0,v) + g(0,v) - g(0,0)  
$$
$$
= \ppp_1g(\xi_2,v)u + \ppp_2g(0,\eta_2)v     \eqno{(6.8)}
$$
where $\xi_2 = \xi_2(u(x,t))$ is a number between $u(x,t)$ and $0$, and
$\eta_2 = \eta_2(v(x,t))$ is a number between $v(x,t)$ and $0$.

Then, the problem (6.1) can be rewritten in the form
$$
\left\{ \begin{array}{rl}
& \pppa (u-a) = \Delta u + p_1(x,t)u + f(0,v(x,t)) \ \ \mbox{in $\OOO\times (0,T)$}, \cr\\
& \pppa (v-b) = \Delta v + q_2(x,t)u + q_3(x,t)v
\ \ \mbox{in $\OOO\times (0,T)$}, \cr\\
& \ppp_{\nu}u = \ppp_{\nu}v = 0 \quad \mbox{on $\ppp\OOO\times (0,T)$},
\end{array}\right.
$$
where we set 
$$
q_2(x,t) := (\ppp_1g)(\xi_2(u(x,t)),v(x,t)), \quad 
q_3(x,t) := (\ppp_2g)(0, \eta_2(v(x,t)).
$$
By the assumptions made for Case 2, the inclusions $q_2, q_3 \in L^{\infty}(\OOO\times (0,T))$
and the inequalities $q_2 (x,t) \ge 0$, $f(0,v(x,t)) \ge 0$ for $(x,t) \in \OOO \times (0,T)$ hold valied. 
Therefore, we can apply Theorem 5.1 that yields non-negativity of the solutions to the problem (6.1):  $u, v \ge 0$ in $\OOO\times (0,T)$.
\\
{\bf Case 3:}
$$
\ppp_2f(\xi,\eta) \ge 0, \quad f(0,0) = 0, \quad g(\xi,0) \ge 0 \quad
\mbox{for $\xi, \eta \in \R$.}      
$$
Then (6.7) holds valid as well as the  representation
$$
f(u(x,t), v(x,t)) = f(u,v) - f(u,0) + f(u,0) - f(0,0)
$$
$$
= \ppp_2f(u,\eta_3)v + \ppp_1f(\xi_3,0)u,     \eqno{(6.9)}
$$
where $\eta_3 = \eta_3(v(x,t))$ is a number between $v(x,t)$ and $0$, and
$\xi_3 = \xi_3(u(x,t))$ is a number between $u(x,t)$ and $0$. Then the arguments as in Case 2 can be applied and we again arrive at non-negativity of the solutions to the problem (6.1). 
\\
{\bf Case 4:}
$$
\ppp_2f(\xi,\eta) \ge 0, \quad f(0,0) = 0, \quad \ppp_1g(\xi,\eta) \ge 0,
\quad g(0,0) = 0 \quad
\mbox{for $\xi, \eta \in \R$.}
$$

Then the representations  (6.8) and (6.9) hold valid and the problem (6.1) can be rewritten in the form
$$
\left\{ \begin{array}{rl}
& \pppa (u-a) = \Delta u + p_2(x,t)u + p_3(x,t)v \ \ \mbox{in $\OOO\times (0,T)$}, \cr\\
& \pppa (v-b) = \Delta v + q_2(x,t)u + q_3(x,t)v
\ \ \mbox{in $\OOO\times (0,T)$}, \cr\\
& \ppp_{\nu}u = \ppp_{\nu}v = 0 \quad \mbox{on $\ppp\OOO\times (0,T)$},
\end{array}\right.
$$
where
$$
p_2(x,t) := (\ppp_1f)(\xi_3(u(x,t)),0), \ 
p_3(x,t) := (\ppp_2f)(u(x,t), \eta_3(v(x,t)).
$$
The assumptions made in Case 4 allow application of Theorem 5.1 to the problem (6.1) in the form presented above that leads to non-negativity of its solutions:  $u, v \ge 0$ in $\OOO\times (0,T)$.

The proof of Theorem 6.1 is completed.
\end{proof} $\square$
\begin{remark}
The semilinear terms $f,\, g$ satisfying the conditions (6.3) and (6.4) often appear in formulations of the mathematical models  of the so called competitive systems, see, e.g., the conditions (4.34) on p. 191 of \cite{Ya} in the form $f(0,\eta) = g(\xi,0) = 0$, which is
contained in our conditions (6.3) and (6.4).  Other examples of this kind can be found 
in \cite{Pi2}. 

We also mention a simple procedure for construction of the semilinear terms satisfying the conditions  of type (6.3) or (6.4). Let $h_k \in C^1(\R^2)$, $1\le k \le 5$.  The function $f$ in form
$$
f(\xi,\eta) = h_1(\xi)h_2(\eta) + h_3(\eta), \ \ 
\mbox{where $h_1(0) = 0$ and $h_3(\eta) \ge 0$ for $\eta\in \R$}
$$
satisfies the first condition in (6.3) and the function $f$ in form
$$
f(\xi,\eta) = h_4(\xi)h_5(\eta), \ \ 
\mbox{where $h_4(0)h_5(0) = 0$ and $h_4(\xi)h_5'(\eta) \ge 0$ for 
$\xi, \eta\in \R$}
$$
satisfies the second condition in (6.3).
\end{remark}

\begin{remark}
The non-negativity of the solutions to the systems of the semilinear diffusion equations plays important role for studying the local existence of solutions in time.  In the case of systems of the parabolic PDEs ($\alpha=1$ in (6.1)), many works are devoted to this subject, see, e.g., 
 \cite{Pi1}, \cite{Pi2},  \cite{QS}, and the
references therein. However, the case of $0<\alpha<1$ has been treated  only in a few works, see, e.g.,  the recent paper \cite{FY}.
\end{remark}

\section{Concluding remarks and directions for further research}
\label{sec_final}

\setcounter{section}{5}
\setcounter{equation}{0}

In this paper, the 
solutions $u\in L^2(0,T;H^2(\OOO))$ to the initial-boundary value problems for 
the semilinear time-fractional diffusion 
equations under the initial condition 
$u-a \in \HH(0,T;L^2(\OOO))$ have been constructed and used for derivation of 
the comparison principles and for justification and applications of the 
monotonicity method. 
However, our main results remain valid for the solutions from  larger spaces 
of functions, for instance, for the mild solutions to the semilinear 
equation \eqref{(2.3)} which  
are defined as solutions to the integral equation \eqref{(7.3)}.

In order to prove the comparison principles for the  mild solutions, 
the pointwise arguments presented in Lemma 3 from \cite{LYP1} and thus 
the approximations of the
solutions by smooth solutions $u\in C([0,T];C^2(\ooo{\OOO}))$
satisfying the inclusion $t^{1-\alpha}\ppp_tu \in C([0,T];C(\ooo{\OOO}))$ 
should be applied. In other words, one first 
proves the comparison principles
for the smooth solutions and then extends them to the desired class of 
solutions. 

In our discussions, we restricted ourselves to the case  of the homogeneous 
boundary conditions.  However, in the definitions
\eqref{(3.10)} and \eqref{(3.11)} of the upper and lower solutions, 
it is natural to  replace 
the homogeneous boundary conditions by the inequalities 
$\NUNU \UPP + \sigma \UPP \ge 0$ and $\NUNU \LOP + \sigma\LOP \le 0$
on $\ppp\OOO \times (0,T)$.  To this end, the unique existence results for the 
initial-boundary value problems with non-homogeneous boundary conditions are 
needed, which are not sufficiently studied. One of the possibilities to attack 
such problems would be to consider the conjugate problems with the homogeneous 
boundary conditions. In this case, the conjugate equations will contain the time-fractional Riemann-Liouville 
derivatives instead of the Caputo derivative. The maximum principle and the non-negativity of solutions to the 
time-fractional diffusion equations with the Riemann-Liouville derivatives 
have been studied in \cite{Al1,Al2,Al3}. However, to the best of our knowledge,
no results regarding existence and uniqueness of solutions to the  
time-fractional diffusion equations with the Riemann-Liouville derivatives and 
with the Neumann or Robin boundary conditions as well as comparison principles 
for solutions to such problems have been yet published. These problems will 
be considered elsewhere.

For the sake of simplicity of arguments, in this paper, we restricted 
ourselves to the case of the spatial dimensions 
$d=1, 2, 3$.  The case $d\ge 4$ can be treated similarly, but some stronger 
regularity conditions
for the coefficients of the differential operators and more involved estimates 
are needed to ensure validity of the results presented in this paper for the 
case $d=1, 2, 3$.




\section*{Acknowledgements}

The work was supported by
Grant-in-Aid for Scientific Research (A) 20H00117 and Grant-in-Aid for 
Challenging Research (Pioneering)
21K18142 of Japan Society for the Promotion of Science.






\end{document}